\documentclass[12pt,reqno]{amsart}

\usepackage{amssymb}

\input psfig.sty

\catcode`\@=11
\newskip\Einheit \Einheit=0.6cm \unitlength=0.6cm
\newcount\xcoord \newcount\ycoord
\newdimen\xdim \newdimen\ydim \newdimen\PfadD@cke \newdimen\Pfadd@cke
\def\PfadDicke#1{\PfadD@cke#1 \divide\PfadD@cke by2 \Pfadd@cke\PfadD@cke \multiply\PfadD@cke by2}
\long\def\LOOP#1\REPEAT{\def\BODY{#1}\ITERATE}
\def\ITERATE{\BODY \let\next\ITERATE \else\let\next\relax\fi \next}
\let\REPEAT=\fi
\def\Punkt{\hbox{\raise-2pt\hbox to0pt{\hss\scriptsize$\bullet$\hss}}}
\def\DuennPunkt(#1,#2){\unskip
  \raise#2 \Einheit\hbox to0pt{\hskip#1 \Einheit
          \raise-2.5pt\hbox to0pt{\hss\normalsize$\bullet$\hss}\hss}}
\def\NormalPunkt(#1,#2){\unskip
  \raise#2 \Einheit\hbox to0pt{\hskip#1 \Einheit
          \raise-3pt\hbox to0pt{\hss\large$\bullet$\hss}\hss}}
\def\DickPunkt(#1,#2){\unskip
  \raise#2 \Einheit\hbox to0pt{\hskip#1 \Einheit
          \raise-4pt\hbox to0pt{\hss\Large$\bullet$\hss}\hss}}
\def\Kreis(#1,#2){\unskip
  \raise#2 \Einheit\hbox to0pt{\hskip#1 \Einheit
          \raise-4pt\hbox to0pt{\hss\Large$\circ$\hss}\hss}}

\def\Diagonale(#1,#2)#3{\unskip\leavevmode
  \xcoord#1\relax \ycoord#2\relax
      \raise\ycoord \Einheit\hbox to0pt{\hskip\xcoord \Einheit
         \line(1,1){#3}\hss}}
\def\AntiDiagonale(#1,#2)#3{\unskip\leavevmode
  \xcoord#1\relax \ycoord#2\relax 
      \raise\ycoord \Einheit\hbox to0pt{\hskip\xcoord \Einheit
         \line(1,-1){#3}\hss}}
\def\Pfad(#1,#2),#3\endPfad{\unskip\leavevmode
  \xcoord#1 \ycoord#2 \ZeichnePfad#3\endPfad}
\def\ZeichnePfad#1{\ifx#1\endPfad\let\next\relax
  \else\let\next\ZeichnePfad
    \ifnum#1=1
      \raise\ycoord \Einheit\hbox to0pt{\hskip\xcoord \Einheit
         \vrule height\Pfadd@cke width1 \Einheit depth\Pfadd@cke\hss}%
      \advance\xcoord by 1
    \else\ifnum#1=2
      \raise\ycoord \Einheit\hbox to0pt{\hskip\xcoord \Einheit
        \hbox{\hskip-\PfadD@cke\vrule height1 \Einheit width\PfadD@cke depth0pt}\hss}%
      \advance\ycoord by 1
    \else\ifnum#1=3
      \raise\ycoord \Einheit\hbox to0pt{\hskip\xcoord \Einheit
         \line(1,1){1}\hss}
      \advance\xcoord by 1
      \advance\ycoord by 1
    \else\ifnum#1=4
      \raise\ycoord \Einheit\hbox to0pt{\hskip\xcoord \Einheit
         \line(1,-1){1}\hss}
      \advance\xcoord by 1
      \advance\ycoord by -1
    \fi\fi\fi\fi
  \fi\next}
\def\hSSchritt{\leavevmode\raise-.4pt\hbox to0pt{\hss.\hss}\hskip.2\Einheit
  \raise-.4pt\hbox to0pt{\hss.\hss}\hskip.2\Einheit
  \raise-.4pt\hbox to0pt{\hss.\hss}\hskip.2\Einheit
  \raise-.4pt\hbox to0pt{\hss.\hss}\hskip.2\Einheit
  \raise-.4pt\hbox to0pt{\hss.\hss}\hskip.2\Einheit}
\def\vSSchritt{\vbox{\baselineskip.2\Einheit\lineskiplimit0pt
\hbox{.}\hbox{.}\hbox{.}\hbox{.}\hbox{.}}}
\def\DSSchritt{\leavevmode\raise-.4pt\hbox to0pt{%
  \hbox to0pt{\hss.\hss}\hskip.2\Einheit
  \raise.2\Einheit\hbox to0pt{\hss.\hss}\hskip.2\Einheit
  \raise.4\Einheit\hbox to0pt{\hss.\hss}\hskip.2\Einheit
  \raise.6\Einheit\hbox to0pt{\hss.\hss}\hskip.2\Einheit
  \raise.8\Einheit\hbox to0pt{\hss.\hss}\hss}}
\def\dSSchritt{\leavevmode\raise-.4pt\hbox to0pt{%
  \hbox to0pt{\hss.\hss}\hskip.2\Einheit
  \raise-.2\Einheit\hbox to0pt{\hss.\hss}\hskip.2\Einheit
  \raise-.4\Einheit\hbox to0pt{\hss.\hss}\hskip.2\Einheit
  \raise-.6\Einheit\hbox to0pt{\hss.\hss}\hskip.2\Einheit
  \raise-.8\Einheit\hbox to0pt{\hss.\hss}\hss}}
\def\SPfad(#1,#2),#3\endSPfad{\unskip\leavevmode
  \xcoord#1 \ycoord#2 \ZeichneSPfad#3\endSPfad}
\def\ZeichneSPfad#1{\ifx#1\endSPfad\let\next\relax
  \else\let\next\ZeichneSPfad
    \ifnum#1=1
      \raise\ycoord \Einheit\hbox to0pt{\hskip\xcoord \Einheit
         \hSSchritt\hss}%
      \advance\xcoord by 1
    \else\ifnum#1=2
      \raise\ycoord \Einheit\hbox to0pt{\hskip\xcoord \Einheit
        \hbox{\hskip-2pt \vSSchritt}\hss}%
      \advance\ycoord by 1
    \else\ifnum#1=3
      \raise\ycoord \Einheit\hbox to0pt{\hskip\xcoord \Einheit
         \DSSchritt\hss}
      \advance\xcoord by 1
      \advance\ycoord by 1
    \else\ifnum#1=4
      \raise\ycoord \Einheit\hbox to0pt{\hskip\xcoord \Einheit
         \dSSchritt\hss}
      \advance\xcoord by 1
      \advance\ycoord by -1
    \fi\fi\fi\fi
  \fi\next}
\def\Koordinatenachsen(#1,#2){\unskip
 \hbox to0pt{\hskip-.5pt\vrule height#2 \Einheit width.5pt depth1 \Einheit}%
 \hbox to0pt{\hskip-1 \Einheit \xcoord#1 \advance\xcoord by1
    \vrule height0.25pt width\xcoord \Einheit depth0.25pt\hss}}
\def\Koordinatenachsen(#1,#2)(#3,#4){\unskip
 \hbox to0pt{\hskip-.5pt \ycoord-#4 \advance\ycoord by1
    \vrule height#2 \Einheit width.5pt depth\ycoord \Einheit}%
 \hbox to0pt{\hskip-1 \Einheit \hskip#3\Einheit 
    \xcoord#1 \advance\xcoord by1 \advance\xcoord by-#3 
    \vrule height0.25pt width\xcoord \Einheit depth0.25pt\hss}}
\def\Gitter(#1,#2){\unskip \xcoord0 \ycoord0 \leavevmode
  \LOOP\ifnum\ycoord<#2
    \loop\ifnum\xcoord<#1
      \raise\ycoord \Einheit\hbox to0pt{\hskip\xcoord \Einheit\Punkt\hss}%
      \advance\xcoord by1
    \repeat
    \xcoord0
    \advance\ycoord by1
  \REPEAT}
\def\Gitter(#1,#2)(#3,#4){\unskip \xcoord#3 \ycoord#4 \leavevmode
  \LOOP\ifnum\ycoord<#2
    \loop\ifnum\xcoord<#1
      \raise\ycoord \Einheit\hbox to0pt{\hskip\xcoord \Einheit\Punkt\hss}%
      \advance\xcoord by1
    \repeat
    \xcoord#3
    \advance\ycoord by1
  \REPEAT}
\def\Label#1#2(#3,#4){\unskip \xdim#3 \Einheit \ydim#4 \Einheit
  \def\lo{\advance\xdim by-.5 \Einheit \advance\ydim by.5 \Einheit}%
  \def\llo{\advance\xdim by-.25cm \advance\ydim by.5 \Einheit}%
  \def\loo{\advance\xdim by-.5 \Einheit \advance\ydim by.25cm}%
  \def\o{\advance\ydim by.25cm}%
  \def\ro{\advance\xdim by.5 \Einheit \advance\ydim by.5 \Einheit}%
  \def\rro{\advance\xdim by.25cm \advance\ydim by.5 \Einheit}%
  \def\roo{\advance\xdim by.5 \Einheit \advance\ydim by.25cm}%
  \def\l{\advance\xdim by-.30cm}%
  \def\r{\advance\xdim by.30cm}%
  \def\lu{\advance\xdim by-.5 \Einheit \advance\ydim by-.6 \Einheit}%
  \def\llu{\advance\xdim by-.25cm \advance\ydim by-.6 \Einheit}%
  \def\luu{\advance\xdim by-.5 \Einheit \advance\ydim by-.30cm}%
  \def\u{\advance\ydim by-.30cm}%
  \def\ru{\advance\xdim by.5 \Einheit \advance\ydim by-.6 \Einheit}%
  \def\rru{\advance\xdim by.25cm \advance\ydim by-.6 \Einheit}%
  \def\ruu{\advance\xdim by.5 \Einheit \advance\ydim by-.30cm}%
  #1\raise\ydim\hbox to0pt{\hskip\xdim
     \vbox to0pt{\vss\hbox to0pt{\hss$#2$\hss}\vss}\hss}%
}
\catcode`\@=12

\def\al{\alpha}

\def\de{\delta}
\def\ep{\varepsilon}
\def\ze{\zeta}

\def\th{\theta}

\def\la{\lambda}

\def\ph{\varphi}

\def\Ga{\Gamma}

\def\Th{\Theta}

\def\Z{{\mathbb Z}}

\def\today{\ifcase\month\or
 January\or February\or March\or April\or May\or June\or
 July\or August\or September\or October\or November\or December\fi
 \space\number\day, \number\year}

\def\({\left(}
\def\){\right)}
\def\[{\left[}
\def\]{\right]}

\def\Tr{\operatorname{Tr}}

\def\Re{\operatorname{Re}}

\def\3{\ss}

\numberwithin{equation}{section}

\newcounter{saveeqn}
\newcommand{\alphaeqn}{\setcounter{saveeqn}{\value{equation}}%
\setcounter{equation}{0}%
\global\def\theequation{\mbox{\thesection.\arabic{saveeqn}\alph{equation}}}}
\newcommand{\reseteqn}{\setcounter{equation}{\value{saveeqn}}%
\global\def\theequation{\thesection.\arabic{equation}}}

\newtheorem{theorem}{Theorem}
\newtheorem{corollary}[theorem]{Corollary}
\newtheorem{lemma}[theorem]{Lemma} 
\newtheorem*{remark}{Remark}

\def\v#1{{\vert #1\vert}}
\def\Tr{\operatorname{Tr}}
\def\po#1#2{(#1)_#2}
\def\fl#1{\left\lfloor#1\right\rfloor}

\def\bX{{\bar X}}

\begin{document}

\author
{Christian Krattenthaler}
\author{Paul B. Slater}
\address{\hskip-\parindent
Christian Krattenthaler
Institut f\"ur Mathematik\\
Universit\"at Wien\\
Strudlhofgasse 4\\
A-1090 Vienna, Austria}
\email{kratt@pap.univie.ac.at \hfil\break\null\hskip\parindent  {\it WWW:}
http://radon.mat.univie.ac.at/People/kratt}
\address{\hskip-\parindent
Paul B. Slater\\
Community and Organization Research Institute\\
University of California\\
Santa Barbara\\
CA 93106-2150}
\email{slater@itp.ucsb.edu}
\title[Asymptotic Redundancies]{Asymptotic Redundancies for Universal\\
Quantum Coding}
\thanks{Krattenthaler's research was supported in part by
MSRI, through NSF grant DMS-9022140.}

\begin{abstract}
Clarke and Barron have recently shown that the Jeffreys' invariant
prior of Bayesian theory yields the common asymptotic
(minimax and maximin) redundancy of universal data compression
in a parametric setting. We seek  a possible analogue
of this result for the two-level {\it quantum} systems.
We restrict our considerations to
 prior
probability distributions belonging
to a certain
 one-parameter family, $q(u)$, $-\infty < u < 1$.
Within this setting, we 
are able to compute exact redundancy formulas, for which we find
 the asymptotic limits.
We compare our
quantum asymptotic redundancy formulas to those derived by
naively applying the classical counterparts
 of Clarke and Barron, and find certain common features.
Our results are based on formulas we obtain for the
eigenvalues and eigenvectors of  $2^n \times 2^n$ (Bayesian 
density) matrices, $\zeta_{n}(u)$.
These matrices are the weighted averages (with respect to $q(u)$)
of all possible tensor products of $n$ identical $2 \times 2$ density matrices,
representing the two-level quantum systems.
We propose a form of {\it universal} coding 
for the situation in which
 the density matrix describing an ensemble of quantum signal states
is unknown. A sequence of $n$ signals would be
projected onto the dominant eigenspaces of $\ze_n(u)$.
\end{abstract}


\maketitle

\section{Introduction}\label{s1}
A theorem has recently been proven \cite{jo,sch} (cf. \cite{ben,cleve,lo}),
in the context of quantum information theory \cite{ben,per}, that is
analogous to the noiseless coding theorem of
classical information theory. 
In the quantum result, the von Neumann entropy \cite{oh,wehrl},
\begin{equation} \label{eq:1}
S(\rho) = -\Tr \rho \log \rho
\end{equation}
(equalling the Shannon entropy of the probability
distribution formed by the eigenvalues of $\rho$) of the density matrix,
\begin{equation} \label{eq:2}
\rho = \sum_{a} p(a) \pi_{a},
\end{equation}
describing an ensemble of {\it pure} quantum signal states,
is equal to $\log 2 \approx .693147$ times the number of quantum bits
(``qubits'') --- that is, the number of two-dimensional Hilbert spaces
--- necessary to represent the signal faithfully.
(Although the binary logarithm is usually used in
the quantum coding literature, we employ the natural logarithm
throughout this paper, chiefly to facilitate comparisons of
our results with those of Clarke and Barron \cite{cl3,cl1,cl2}. 
$p(a)$ is the probability of the message $a$ from a particular
source coded into a ``signal state'' --- having a state
vector denoted by the ket $|a_{M}\rangle$ --- of a quantum system $M$.
The density matrices $\pi_{a}$ are the projections
$\pi_{a} = |a_{M}\rangle \langle a_{M}|$, with $\langle a_{M}|$ being a bra
in the dual Hilbert space.)

 The proof
of the quantum coding theorem is based on the existence
of a ``typical subspace'' $\Lambda$ of the
$2^n$-dimensional Hilbert space of
$n$ qubits, which has the property that, with
high probability, a sample of $n$ qubits has almost
unit projection onto $\Lambda$. Since it has
been shown that the dimension of $\Lambda$ is
$e^{n S(\rho)}$, the operation that the
data compressor (a unitary transformation mapping
$n$-qubit strings to $n$-qubit strings) should
perform involves ``transposing'' the subspace $\Lambda$ into
the Hilbert space of a smaller block of
${n S(\rho)} /{.693147}$ qubits \cite{cleve}.
(Lo \cite{lo} has generalized this
work for an ensemble of {\it mixed} quantum signal states.)

In this study we dispense with the assumption that
{\it a priori} information (other
than its dimensionality) is available regarding $\rho$.
Somewhat similarly motivated, Calderbank and Shor \cite{cald} modified
the definition of {\it fidelity} --- a measure of the success of
transmission of quantum states --- because ``previous papers discuss
channels that transmit some distribution of states given {\it a priori},
whereas we want our channel to faithfully transmit any pure input
state''. They took as their measure, the fidelity for the pure state
transmitted {\it least} faithfully.

Proceeding in a noninformative
Bayesian framework \cite{bern,sl1,sl2,sl3}, we seek to extend to the two-level quantum systems,
recent results of Clarke and Barron \cite{cl3,cl1,cl2}
giving various forms of the asymptotic redundancy of universal
data compression for parameterized families of
probability distributions.
``The redundancy is the excess of the [coding] cost over
the entropy. The goal of data compression is to diminish
redundancy'' 
(\cite{kir}, reviewed in \cite{csis}).
``The idea of universal coding, suggested by Kolmogorov,
is to construct a code for data sequences such that asymptotically,
as the length of the sequence increases, the mean per symbol code
length would approach the entropy of whatever process in a family
has generated the data'' \cite{ri}.
For an extensive commentary on the
results of Clarke and Barron, see \cite{ri}. Also see \cite{cl4}, for
some recent related research, as well as a
discussion of various rationales
that have been employed for using the (classical) Jeffreys'
 prior --- a possible
quantum counterpart of which will be of interest here --- for
Bayesian purposes, cf. \cite{kass}.
Let us also bring to the attention
of the reader that in a brief review of \cite{cl1}, the noted
statistician, I. J. Good, commented that Clarke and Barron `` have presumably
overlooked the reviewer's work'' and cited, in this regard
 \cite{good2,good3}.
(It should be noted that in these papers, Good uses a more general
objective function --- a two-parameter utility --- than the relative
entropy, chosen by Clarke and Barron over
alternative measures \cite[p.~454]{cl3}. Good does conclude that Jeffreys'
invariant prior is the minimax, that is, the least favorable, prior when
the utility is the ``weight of the evidence'' in the sense of C. S.
Pierce, that is, the relative entropy.)

Clarke and Barron \cite{cl3,cl1,cl2} found
the asymptotic redundancy to be given by
\begin{equation} \label{eq:4}
{\frac d  2} {\log {\frac n  {2 \pi e}}}+{\frac 1  2}\log \det
I(\theta)-\log w(\theta) +o(1). 
\end{equation}
Here, $\theta$ is a $d$-dimensional vector of variables
parameterizing a family (manifold) of probability distributions. $I(\theta)$
is the $d \times d$ Fisher information matrix --- the
negative of the expected value of the Hessian
of the logarithm of the density function ---  
and $w(\theta)$ is the prior 
density.
The asymptotic {\it minimax\/}
redundancy was shown to be \cite{cl1,cl2}
\begin{equation} \label{eq:3}
{\frac d  2} {\log {\frac n  {2 \pi e}}} +\log {\int_{K} \sqrt{\det
I(\theta)}} \,
d \theta + o(1),
\end{equation}
where $K$ is a compact
set in the interior of the
domain of the parameters.

In this investigation, instead of probability 
densities
as in \cite{cl3,cl1,cl2}, we employ
density matrices (nonnegative definite Hermitian matrices
of unit trace) and instead of the classical form of the relative entropy
(the Kullback--Leibler information measure),
its quantum counterpart \cite{oh,wehrl} (cf. \cite{pe4}),
\begin{equation} \label{eq:5}
 S(\rho_{1},\rho_{2}) =
 \Tr \rho_{1} (\log \rho_{1} -  \log \rho_{2}),
\end{equation}
that is, the relative entropy of the density matrix $\rho_{1}$ with respect
to $\rho_{2}$.

The three-dimensional convex set of $2\times 2$ density matrices
that will be the focus of our study has members representable in the form,
\begin{equation} \label{eq:6}
\rho = {\frac 1  2} \begin{pmatrix}1+z&x-iy\\
                   x+iy&1-z\end{pmatrix}\quad .
\end{equation}
Such matrices correspond, in a one-to-one fashion, to the standard 
(complex) two-level quantum systems --- notably, those of
spin-$1/2$ (electrons, protons,\dots)
and massless spin-$1$ particles (photons).
 (If we set $x=y=0$ in (\ref{eq:6}), we recover a
classical binomial distribution, with the probability of ``success'', say, being $(1+z)/2$ and of ``failure'', $(1-z)/2$. Setting either
$x$ or $y$ to zero, puts us in the
framework of real --- as opposed to complex --- quantum
mechanics.) The points $(x,y,z)$ must lie within
the unit ball (``Bloch sphere'' \cite{br}), $x^2 +y^2+z^2 \leq 1$, due to the requirement for $\rho$ of
nonnegative eigenvalues. 
(The points on the bounding spherical surface, $x^2+y^2+z^2=1$,
corresponding to the pure states, will be shown 
to exhibit nongeneric behavior, 
see (\ref{a5}) and the respective comments in Sec.~\ref{s3}
(cf. \cite{fuj}).)
We have, for (\ref{eq:6}), using spherical coordinates $(r,\vartheta,\phi)$,
 so that 
$r=(x^2+y^2+z^2)^{1/2}$,
\begin{equation} \label{eq:7}
S(\rho)=-{\frac{(1-r)}  2}{ \log {\frac{(1-r)}  2}} -{\frac {(1+r) }
2}{ \log {\frac{(1+r)}   2}} .
\end{equation}
A composite system of $n$ identical independent
(unentangled)
two-level quantum systems is represented by the $2^n\times 2^n$ density matrix
$\overset n  \otimes \rho$ --- possessing a von Neumann entropy
$n S(\rho)$ \cite{oh,wehrl}.
(In noncommutative probability theory, independence can be based
on free products instead of tensor products \cite{sza}.
Along with the real and complex forms of quantum mechanics, a quaternionic
version exists \cite{fivel}, for which the [presumed] quantum
Jeffreys' prior has been found
for the two-level systems --- corresponding
to the {\it five}-dimensional
unit ball/``Bloch sphere'' \cite{sl1}. However, the definition of a tensor
product is somewhat problematical in this context \cite{adl,del}.)

In \cite{sl1}  it was argued that the quantum Fisher information
matrix (requiring --- due to noncommutativity ---
the computation of symmetric logarithmic
derivatives \cite{pe1}) for the density matrices (\ref{eq:6}) should
be taken to be of the form
\begin{equation} \label{eq:8}
I(\theta) =
{\frac1 {(1-x^2-y^2-z^2)}} \begin{pmatrix}1-y^2-z^2&xy&xz\\
                                 xy&1-x^2-z^2&yz\\
                                 xz&yz&1-x^2-y^2\end{pmatrix}\quad.
\end{equation}
The quantum counterpart of the
 Jeffreys' prior was, then, taken to be the normalized form
(dividing by $\pi^2$)
of the square root of the determinant of (\ref{eq:8}), that is,
\begin{equation} \label{eq:9}
 (1-x^2-y^2-z^2)^{-1/2}/{\pi}^2.
\end{equation}
Analogously, the {\it classical} Jeffreys' prior is proportional to the
square root of the determinant of the {\it classical} Fisher
information matrix \cite{bern}.

On the basis of the result of Clarke and Barron \cite{cl1,cl2} that
the Jeffreys' prior yields the asymptotic common (minimax
and maximin) redundancy (that is, the least favorable and
reference priors are the same),
it was
conjectured
\cite{sl4} that its assumed quantum counterpart (\ref{eq:9}) would
have similar properties,
as well.
(The Jeffreys' prior has been ``shown to be a minimax solution in 
a --- two person --- zero sum game, where the statistician
chooses the `non-informative' prior
and nature chooses the `true' prior'' \cite{bern,kash}.
Quantum mechanics itself has been asserted to arise from
a Fisher-information transfer zero sum game \cite{fried}.)
To examine this possibility,
(\ref{eq:9}) was embedded as a specific member ($u=.5$) of a
one-parameter family of 
spherically-symmetric/unitarily-invariant probability 
densities,
\begin{equation} \label{eq:10}
q(u)=\frac {\Gamma(5/2-u)} {\pi^{3/2}\,\Gamma(1-u)\,(1-x^2-y^2-z^2)^u},\quad -\infty<u<1 .
\end{equation}
(Under unitary transformations of $\rho$, the assigned probability is 
invariant.)
For $u=0$, we obtain a uniform distribution over the unit ball.
(This has been used as a prior over the two-level quantum systems,
at least, in one study \cite{larson}.) For $u \rightarrow 1$,
the uniform distribution over the spherical boundary
(the locus of the pure states) is approached.
(This is often employed as a prior, for example \cite{jones,larson,massar}.)
For $u \rightarrow -\infty$, a Dirac distribution concentrated at
the origin (corresponding to the fully mixed state) is approached.

Embeddings of (\ref{eq:9}) in other (possibly, multiparameter)
families are, of course, possible and may be pursued
in further research. Ideally, we would aspire to formally
demonstrate --- if it is, in fact, so --- that
(\ref{eq:9}) can be uniquely characterized {\it vis-\`a-vis}
{\it all} other possible probability distributions over the unit ball.
Due to the present lack of any such fully rigorous treatment,
analogous to that of Clarke and Barron,
we rely upon an
exploratory heuristic computational strategy.
This involves averaging $\overset n\otimes \rho$ with respect to
$q(u)$. Doing so yields a one-parameter family 
of $2^n \times 2^n$ Bayesian density matrices
(Bayes codes or estimators \cite{cl2,cl3,mat}), $\zeta_{n}(u)$,
$-\infty < u < 1$,
exhibiting highly interesting properties.

We explicitly find (in Sec.~\ref{s2}) 
the eigenvalues and eigenvectors of the matrices
$\ze_n(u)$ and determine the relative entropy (\ref{eq:5}) of
$\overset n\otimes \rho$ with respect to $\ze_n(u)$. 
We do this by using identities for hypergeometric series and some
combinatorics. (It is also possible to obtain some of our results by
making use of representation theory of $SU(2)$. An even more general
result was derived by combining these two approaches. We comment on
this issue at the end of Sec.~\ref{s3}.)

The matrices $\ze_n(u)$ should prove useful for the {\it universal}
version of Schumacher data compression 
\cite{ben,cleve,jo,sch} 
by projecting blocks of $n$ signals (qubits) onto those
``typical'' subspaces 
of $2^n$-dimensional Hilbert space corresponding to as many of
the dominant eigenvalues of $\zeta_{n}(u)$ as it takes to exceed a
sum $1- \epsilon$.
(This can be accomplished by a unitary transformation, the inverse
of which would be used in the decoding step \cite{ben}.
In the corresponding nonuniversal quantum coding context,
the projection onto the dominant eigenvalues of $\overset n\otimes \rho$ yields
fidelity greater than $1-2 \epsilon$ \cite{jo}
and distortion less than $2 \epsilon$ \cite{lo}, cf. \cite{barn}.)
For all $u$, the leading 
one of the $\fl {\frac {n} {2}} + 1$ 
distinct eigenvalues has multiplicity $n+1$, and belongs to the
($n+1$)-dimensional (Bose--Einstein) symmetric subspace \cite{Ba1}.
(Projection onto the symmetric subspace has been proposed as a
method for stabilizing quantum computations, including quantum
state storage \cite{bar}.)
For $u=1/2$, the leading eigenvalue can be obtained by dividing
the $n+1$-st Catalan number 
--- that is, 
$\frac {1} {n+2}\binom {2(n+1)}{n+1}$ --- by $4^n$.
(The Catalan numbers ``are probably the most frequently occurring
combinatorial numbers after the binomial coefficients'' \cite{sloane}.)

Let us (naively) attempt to apply the formulas
of Clarke and Barron \cite{cl1,cl2} --- (\ref{eq:3}) and (\ref{eq:4})
above --- to
the quantum context under investigation here. We do this by
setting $d$ to 3 (the dimensionality of the unit ball --- which we
take as $K$), $\det I(\theta)$ to $(1-x^2-y^2-z^2)^{-1}$ (cf. (\ref{eq:8})),
so that $\int_{K} \sqrt {\det I(\theta)}\,d\th$ is $\pi^2$, and
$w(\theta)$ to $q(u)$. Then, we obtain from the expression for
the asymptotic {\it minimax} redundancy (\ref{eq:3}),
\begin{equation} \label{eq:11}
{\frac 3  2} (\log n -\log 2 -1) + {\frac 1  2} \log \pi + o(1),
\end{equation}
and from the expression for the asymptotic redundancy itself (\ref{eq:4}),
\begin{equation} \label{eq:12}
{\frac 3  2} (\log n - \log 2 -1) -(1-u) \log (1-r^2) + \log \Gamma (1-u)
- \log \Gamma\({\frac 5  2} - u\) + o(1)
\end{equation}
We shall
(in Sec.~\ref{s3})
compare these two formulas,
(\ref{eq:11}) and (\ref{eq:12}), with the results of Sec.~\ref{s2}
and find some striking similarities
and coincidences, particularly associated with the
fully mixed state ($r=0$).
These findings will help to support the working hypothesis
of this study --- that there are meaningful extensions to the
quantum domain of the 
(commutative probabilistic) theorems of Clarke and Barron.
However, we find that although the minimax property of
the Jeffreys' prior appears to carry over, the maximin
property does not strictly, but only in an approximate sense.
In any case, we can not formally rule out the possibility that the
actual global (perhaps common) minimax and maximin are achieved for
probability distributions not belonging to the one-parameter
family $q(u)$.

Let us point out to the reader the
quite recent important work of Petz and Sudar \cite{pe1}.
They demonstrated that in the quantum case --- in
 contrast to the classical situation
in which there is, as originally shown by Chentsov \cite{chen}, essentially only one monotone metric and,
therefore, essentially only one form of the Fisher information --- there
 exists an infinitude of such metrics.
``The monotonicity of the Riemannian metric $g$ is crucial when one likes to
imitate the geometrical approach of [Chentsov]. An infinitesimal statistical
 distance
has to be monotone under stochastic mappings. We note that the monotonicity
of $g$ is a strengthening of the concavity of the von Neumann entropy.
Indeed, positive definiteness of $g$ is equivalent to the strict concavity
of the von Neumann entropy
\dots{} and monotonicity is much more than positivity''
\cite{pe3}.

 The monotone metrics on the space of density matrices are
given \cite{pe1} by the operator monotone functions $f(t):
\mathbb R^+ \rightarrow
\mathbb R^+$, such that
$f(1)=1$ and $f(t) = t f(1/t)$. For the choice $f= (1+t)/2$,
one obtains the minimal metric (of the symmetric logarithmic
derivative), which serves as the basis of our analysis here.
``In accordance with the work of Braunstein and Caves, this seems
to be the canonical metric of parameter estimation theory. However,
expectation values of certain relevant observables are known to lead
to statistical inference theory provided by the maximum entropy
principle or the minimum relative entropy principle when {\it a priori}
information on the state is available. The best prediction is a kind of
generalized Gibbs state. On the manifold of those states, the
differentiation of the entropy functional yields the
 Kubo-Mori/Bogoliubov metric,
which is different from the metric of the symmetric logarithmic
derivative. Therefore, more than one privileged metric shows up
in quantum mechanics. The exact clarification of this point requires
and is worth further studies'' \cite{pe1}.
It remains a possibility, then, that a monotone metric other than the
minimal one (which corresponds to $q(.5)$, that is (\ref{eq:9})) may yield
a  common global asymptotic minimax and maximin
 redundancy, thus, fully paralleling
the classical/nonquantum results of Clarke and Barron
\cite{cl3,cl1,cl2}. We intend to investigate such a possibility,
in particular, for the Kubo-Mori/Bogoliubov metric \cite{pe3,pe1,pe2}.

\section{Analysis of a One-Parameter Family of Bayesian Density Matrices}
\label{s2}
In this section, we implement the analytical approach described in the
Introduction to extending the work of Clarke and Barron \cite{cl1,cl2}
to the realm of quantum mechanics, specifically, the two-level
systems. Such systems are representable by density matrices $\rho$ of the
form (\ref{eq:6}). A composite system of $n$ independent
(unentangled) and identical
two-level quantum systems is, then, represented by the $n$-fold tensor
product $\overset n\otimes \rho$. In
Theorem~\ref{t1} of Sec.~\ref{s2.1},
we average
$\overset n\otimes\rho$ with respect to the one-parameter family of 
probability 
densities $q(u)$ defined in (\ref{eq:10}), obtaining
the Bayesian density matrices $\ze_n(u)$
and formulas for their $2^{2 n}$ entries. Then, in
Theorem~\ref{t2}
of Sec.~\ref{s2.2}, we are able to explicitly determine the $2^n$ eigenvalues
and eigenvectors of $\ze_n(u)$. Using these results, in
Sec.~\ref{s2.3}, we compute the relative entropy of $\overset
n\otimes\rho$ with respect to $\ze_n(u)$. Then, in Sec.~\ref{s2.4}, we
obtain the asymptotics of this relative entropy for $n\to\infty$. 
In Sec.~\ref{s2.5}, we compute the asymptotics of the von Neumann
entropy (see (\ref{eq:1})) of $\ze_n(u)$.
All these results 
will enable us, in Sec.~\ref{s3}, to ascertain to what extent the
results of Clarke and Barron
could be said to carry over to the quantum domain.

\subsection{Entries of the Bayesian density matrices
$\ze_n(u)$}\label{s2.1}
The $n$-fold tensor product $\overset n\otimes \rho$ is a
$2^n\times 2^n$ matrix. To refer to specific rows and columns of 
$\overset n\otimes \rho$, we index them by subsets of
the $n$-element set $\{1,2,\dots,n\}$. We choose to employ this
notation instead of the more familiar use of binary strings, in order 
to have a
more succinct way of writing our formulas. For convenience, 
we will subsequently write $[n]$ for $\{1,2,\dots,n\}$.
Thus, $\overset n\otimes \rho$ can be written in the form
$$\overset n\otimes \rho=\begin{pmatrix} R_{IJ}\end{pmatrix}_{I,J\in
[n]},$$
where
\begin{equation}\label{e2}
R_{IJ}=\frac {1} {2^n}(1+z)^{n_{\in\in}} (1-z)^{n_{\notin\notin}}
(x+iy)^{n_{\notin\in}} (x-iy)^{n_{\in\notin}},
\end{equation}
with $n_{\in\in}$ denoting the number of elements of $[n]$
contained in both $I$ and $J$, $n_{\notin\notin}$ denoting the number
of elements {\it not\/} in both $I$ and $J$, $n_{\notin\in}$ denoting the number
of elements not in $I$ but in $J$, and $n_{\in\notin}$ denoting the number
of elements in $I$ but not in $J$. In symbols,
\begin{align*} n_{\in\in}&=\v{I\cap J},\\ n_{\notin\notin}&=\v{[n]
\backslash (I\cup J)},\\ n_{\notin\in}&=\v{J\backslash I},\\
n_{\in\notin}&=\v{I\backslash J}.
\end{align*}
We consider the average $\ze_n(u)$ of $\overset n\otimes \rho$ with respect
to the probability 
density $q(u)$ defined in (\ref{eq:10})
taken over the unit sphere $\{(x,y,z):x^2+y^2+z^2\le 1\}$.
This average can be described explicitly as follows.
\begin{theorem}\label{t1}
The average $\ze_n(u)$,
$$\int_{x^2+y^2+z^2\le 1}\big(\overset n\otimes \rho\big)\,
q(u)\,dx\,dy\,dz,$$ 
equals the matrix $(Z_{IJ})_{I,J\in[n]}$,
where
\begin{multline} \label{e4}
Z_{IJ}=\de_{n_{\notin\in},n_{\in\notin}}\big(\tfrac {n-n_{\in\in}-
n_{\notin\notin}} 2\big)!\,\\
\times
\frac {1} {2^n}\frac {\Ga\(\frac {5} {2}-u\)\,\Ga\(2+\frac {n} {2}+
\frac{n_{\in\in}}2-\frac {n_{\notin\notin}}2-u\)\, 
\Ga\(2+\frac {n} {2}+\frac {n_{\notin\notin}}2-
\frac {n_{\in\in}}2-u\)} {\Ga\(\frac {5} {2}+\frac {n} {2}-u\)\,
\Ga\(2+\frac {n} {2}-u\)\, \Ga\(2+\frac {n} {2}-\frac {n_{\in\in}}2-
\frac {n_{\notin\notin}}2-u\)}.
\end{multline}
Here, $\de_{i,j}$ denotes the Kronecker delta, $\de_{i,j}=1$ if $i=j$
and $\de_{i,j}=0$ otherwise.
\end{theorem}

\begin{remark} \em It is important for later considerations to observe
that because of the term $\de_{n_{\notin\in},n_{\in\notin}}$ in (\ref{e4})
the entry $Z_{IJ}$ is nonzero if and only if the sets $I$ and $J$
have the same cardinality. If $I$ and $J$ have the same cardinality,
$c$ say,
then $Z_{IJ}$ only depends on $n_{\in\in}$, the number of common
elements of $I$ and $J$, since in this case $n_{\notin\notin}$ is
expressible as $n-2c+n_{\in\in}$.
\end{remark}
\medskip
{\sc Proof of Theorem~\ref{t1}}.
To compute $Z_{IJ}$, we have to compute the integral 
\begin{equation}\label{e5}
\int_{x^2+y^2+z^2\le 1}R_{IJ}\,q(u)\,dx\,dy\,dz.
\end{equation}
For convenience, we treat the case that $n_{\in\in}\ge n_{\notin\notin}$
and $n_{\notin\in}\ge n_{\in\notin}$. The other four cases are
treated similarly.

First, we rewrite the matrix entries $R_{IJ}$,
\begin{align} \notag
\frac {1} {2^n}(1+z)^{n_{\in\in}}& (1-z)^{n_{\notin\notin}}
(x+iy)^{n_{\notin\in}} (x-iy)^{n_{\in\notin}}\\
\notag
=&\frac {1} {2^n}(1-z^2)^{n_{\in\in}}
(1-z)^{n_{\notin\notin}-n_{\in\in}}
(x^2+y^2)^{n_{\notin\in}} (x-iy)^{n_{\in\notin}-n_{\notin\in}}\\
\notag
=&\frac {1} {2^n}\sum _{j,k,l\ge0} ^{}(-1)^{j+k}(-i)^l
\binom {n_{\in\in}}j \binom
{n_{\notin\notin}-n_{\in\in}}k \binom
{n_{\in\notin}-n_{\notin\in}}l \\
&\hskip2cm\cdot z^{2j+k}(x^2+y^2)^{n_{\notin\in}}
x^{n_{\in\notin}-n_{\notin\in}-l} y^l.\label{e6}
\end{align}
Of course, in order to compute the integral (\ref{e5}), we transform the
Cartesian coordinates into polar coordinates,
\begin{align}\notag x&=r\sin\vartheta \cos\ph\\
y&=r\sin\vartheta \sin\ph\notag\\
z&=r\cos\vartheta,\notag\\
0\le \ph\le {}&2\pi,\ 0\le\vartheta\le\pi.\notag
\end{align}
Thus, using (\ref{e6}), the integral (\ref{e5}) is transformed into
\begin{multline} \label{e7}
\frac {1} {2^n}\sum _{j,k,l\ge0} ^{}
\int_0^1\int_0^{\pi}\int_0^{2\pi} (-1)^{j+k}(-i)^l
\binom {n_{\in\in}}j \binom
{n_{\notin\notin}-n_{\in\in}}k \binom
{n_{\in\notin}-n_{\notin\in}}l \\
\hskip2cm\cdot r^{2j+k+n_{\notin\in}+n_{\in\notin}+2}
\(\cos^{2j+k}\vartheta\)\(\sin^{n_{\notin\in}+n_{\in\notin}+1}\vartheta \)\\
\hskip2cm\cdot \(\cos^{n_{\in\notin}-n_{\notin\in}-l}\ph \)\(\sin^l\ph\)
\frac {\Ga(5/2-u)} {\pi^{3/2}\,\Ga(1-u)\,(1-r^2)^u}
\,d\ph\,d\vartheta\,dr.
\end{multline}
To evaluate this triple integral we use the following standard
formulas:
{\refstepcounter{equation}\label{e8}}
\alphaeqn
\begin{align} \int_0^\pi\sin^{2M}\vartheta\,
\cos^{2N}\vartheta\,d\vartheta&=\pi\frac {(2M-1)!!\,(2N-1)!!}
{(2M+2N)!!},\label{e8a}\\
\int_0^\pi\sin^{2M+1}\vartheta\,
\cos^{2N}\vartheta\,d\vartheta&=2\frac {(2M)!!\,(2N-1)!!}
{(2M+2N+1)!!}\ \\
&\text {and}\quad 
\int_0^{2\pi}\sin^{2M+1}\vartheta\,
\cos^{2N}\vartheta\,d\vartheta=0,\label{e8b}\\
\int_0^\pi\sin^{2M}\vartheta\,
\cos^{2N+1}\vartheta\,d\vartheta&=0,\label{e8c}\\
\int_0^\pi\sin^{2M+1}\vartheta\,
\cos^{2N+1}\vartheta\,d\vartheta&=0,\label{e8d}
\end{align}
\reseteqn
for any nonnegative integers $M$ and $N$. Furthermore, we need
the beta integral
\begin{equation}\label{e9}
\int_0^1\frac {r^m} {(1-r^2)^u}\,dr=\frac {\Ga\( \frac {m+1}
{2}\)\,\Ga(1-u)} {2\,\Ga\(\frac {m+3} {2}-u\)}.
\end{equation}

Now we consider the integral over $\ph$ in (\ref{e7}). Using
(\ref{e8b}) and (\ref{e8c}),
we see that each summand in (\ref{e7}) vanishes if $n_{\notin\in}$ has a
parity different from $n_{\in\notin}$. On the other hand, if 
$n_{\notin\in}$ has the same parity as $n_{\in\notin}$, then we can
evaluate the integrals over $\ph$ using (\ref{e8a}) and
(\ref{e8d}). 
Discarding for a
moment the terms independent of $\ph$ and $l$, we have
\begin{align} \sum _{l\ge0} ^{}&\int_0^{2\pi} (-i)^l  \binom
{n_{\in\notin}-n_{\notin\in}}l 
\(\cos^{n_{\in\notin}-n_{\notin\in}-l}\ph \)\(\sin^l\ph\)\,d\ph\notag\\
&=\sum _{l\ge0} ^{}(-1)^l\binom
{n_{\in\notin}-n_{\notin\in}}{2l} 2\pi\,\frac {(2l-1)!!\,(n_{\in\notin}
-n_{\notin\in}-2l-1)!!} {(n_{\in\notin}-n_{\notin\in})!!}\notag\\
&=2\pi\frac {(n_{\in\notin}-n_{\notin\in}-1)!!} 
{(n_{\in\notin}-n_{\notin\in})!!}\sum _{l\ge0} ^{}\binom {(n_{\in\notin}-n_{\notin\in})/2}l
(-1)^l\notag\\
&=2\pi\,\de_{n_{\in\notin},n_{\notin\in}},\notag
\end{align}
the last line being due to the binomial theorem. These 
considerations reduce (\ref{e7}) to
\begin{multline}\notag \de_{n_{\in\notin},n_{\notin\in}}\frac {1} {2^n}
\sum _{j,k\ge0} ^{}
\int_0^1\int_0^{\pi} (-1)^{j+k}
\binom {n_{\in\in}}j \binom
{n_{\notin\notin}-n_{\in\in}}k \\
\hskip2cm\cdot r^{2j+k+2n_{\notin\in}+2}
\(\cos^{2j+k}\vartheta\)\(\sin^{2n_{\notin\in}+1}\vartheta \)
 \frac {2\,\Ga(5/2-u)} {\pi^{1/2}\,\Ga(1-u)\,(1-r^2)^u}
\,d\vartheta\,dr.
\end{multline}
Using (\ref{e8b}), (\ref{e8d}) and (\ref{e9}) this can be further simplified to
\begin{multline}\label{e10} 
\de_{n_{\in\notin},n_{\notin\in}}\frac {1} {2^n}
\sum _{j,k\ge0} ^{}
 (-1)^{j}
\binom {n_{\in\in}}j \binom
{n_{\notin\notin}-n_{\in\in}}{2k}  \frac {2\,(2j+2k-1)!!\,(2n_{\notin\in})!!}
{(2j+2k+2n_{\notin\in}+1)!!}\\
\cdot\frac
{\Ga(j+k+n_{\notin\in}+3/2)\,\Ga(1-u)}
{2\,\Ga(j+k+n_{\notin\in}+5/2-u)}
 \frac {2\,\Ga(5/2-u)} {\pi^{1/2}\,\Ga(1-u)}.
\end{multline}
Next we interchange sums over $j$ and $k$ 
and write the sum over $k$ in terms of the standard
hypergeometric notation
$${}_r F_s\!\left[\begin{matrix} a_1,\dots,a_r\\ b_1,\dots,b_s\end{matrix}; 
z\right]=\sum _{k=0} ^{\infty}\frac {\po{a_1}{k}\cdots\po{a_r}{k}}
{k!\,\po{b_1}{k}\cdots\po{b_s}{k}} z^k\ ,$$
where the shifted factorial
$(a)_k$ is given by $(a)_k:=a(a+1)\cdots(a+k-1)$,
$k\ge1$, $(a)_0:=1$. Thus we can write (\ref{e10}) in the form
\begin{multline}\label{e11} 
\de_{n_{\in\notin},n_{\notin\in}}\frac {1} {2^n}
\sum _{k\ge0} ^{}\binom {n_{\notin\notin}-n_{\in\in}}{2k}
\frac {(2k-1)!!\,n_{\notin\in}!\,\Ga\(\frac {5} {2}-u\)} {2^{k+1}\,\Ga\(\frac
{5} {2}+k+n_{\notin\in}-u\)}\\
\cdot{}_2F_1\!\[\begin{matrix} \frac {1} {2}+k,-n_{\in\in}\\ \frac {5} {2}+k+n_{\notin\in}-u
\end{matrix}; 1\].
\end{multline}
The $_2F_1$ series can be summed by means of Gau\ss' $_2F_1$
summation (see e.g\@. \cite[(1.7.6); Appendix (III.3)]{SlatAC})
\begin{equation}\label{e12}
{} _{2} F _{1} \!\left [ \begin{matrix} { a, b}\\ { c}\end{matrix} ; {\displaystyle
   1}\right ]  = \frac {\Ga(c)\,\Ga(c-a-b)} {\Ga(c-a)\,\Ga(c-b)} ,
\end{equation}
provided the series terminates or $\Re (c-a-b)\ge0$. Applying (\ref{e12}) to
the $_2F_1$ in (\ref{e11}) (observe that it is terminating)
and writing the sum over $k$ as a hypergeometric
series, the
expression (\ref{e11}) becomes
\begin{multline}\notag\de_{n_{\in\notin},n_{\notin\in}}\frac {1} {2^n}
\frac {\Ga(2+n_{\in\in}+n_{\notin\in}-u)\,\Ga\(\frac {5}
{2}-u\)\, n_{\notin\in}!} {\Ga\(\frac {5} {2}+n_{\in\in}+n_{\notin\in}-u\)
\,\Ga(2+n_{\notin\in}-u)}\\
\times {}_2F_1\!\[\begin{matrix} \frac {n_{\in\in}} {2}-\frac {n_{\notin\notin}} 
{2},\frac {1} {2}+\frac {n_{\in\in}} {2}-\frac {n_{\notin\notin}} 
{2}\\\frac {5} {2}+n_{\in\in}+n_{\notin\in}-u\end{matrix}; 1\].
\end{multline}
Another application of (\ref{e12}) gives
\begin{multline}\label{e13}
\de_{n_{\in\notin},n_{\notin\in}}\frac {1} {2^n}\\
\times
\frac {\Ga(2+n_{\in\in}+n_{\notin\in}-u)\,
\Ga(2+n_{\notin\notin}+n_{\notin\in}-u)\,\Ga\(\frac {5}
{2}-u\)\, n_{\notin\in}!} {\Ga\(\frac {5} {2}+\frac {n_{\in\in}} {2}+
\frac {n_{\notin\notin}} {2}+n_{\notin\in}-u\)\,\Ga\( {2}+\frac {n_{\in\in}} {2}+
\frac {n_{\notin\notin}} {2}+n_{\notin\in}-u\)
\,\Ga(2+n_{\notin\in}-u)}.\\
\end{multline}
Trivially, we have $n=n_{\in\in}+n_{\notin\notin}+n_{\notin\in}
+n_{\in\notin}$.
Since (\ref{e13}) vanishes unless $n_{\notin\in}=n_{\in\notin}$, we
can substitute
$(n-n_{\in\in}-n_{\notin\notin})/2$ for $n_{\notin\in}$ in the
arguments of the gamma functions. Thus,
we see that (\ref{e13}) equals (\ref{e4}). This completes the proof of the Theorem.\quad \quad
\qed
\medskip

\subsection{Eigenvalues and eigenvectors of the Bayesian density
matrices $\ze_n(u)$}\label{s2.2}

With the explicit description of the result $\ze_n(u)$ of averaging 
$\overset n\otimes\rho$ with respect to $q(u)$
at our disposal, we now proceed to describe the eigenvalues and eigenspaces
of $\ze_n(u)$. The eigenvalues are given in
Theorem~\ref{t2}. Lemma~\ref{l4} gives a complete set of eigenvectors
of $\ze_n(u)$. The reader should note that, though complete,
this is simply a set of linearly independent eigenvectors and not a
fully orthogonal set.
\begin{theorem}\label{t2}
The eigenvalues of the $2^n\times2^n$ matrix 
$\ze_n(u)$, the entries of which are given by {\em(\ref{e4})},
are
\begin{equation}\label{e14}
\la_d=\frac {1} {2^n}\frac {\Ga\(\frac {5} {2}-u\)\, \Ga(2+n-d-u)\,
\Ga(1+d-u)} {\Ga\(\frac {5} {2}+\frac {n} {2}-u\)\, \Ga(2+\frac {n}
{2}-u)\,\Ga(1-u)}, \quad d=0,1,\dots,\fl{\frac {n} {2}},
\end{equation}
with respective multiplicities
\begin{equation}\label{e15}
\frac {(n-2d+1)^2} {(n+1)}\binom {n+1}d.
\end{equation}
\end{theorem}
The Theorem will follow from a sequence of Lemmas. We state the
Lemmas first, then prove Theorem~\ref{t2} assuming the truth of the Lemmas,
and after that provide proofs of the Lemmas.

In the first Lemma some eigenvectors of the matrix $\ze_n(u)$ are
described. Clearly, since $\ze_n(u)$ is a $2^n\times2^n$ matrix, the
eigenvectors are in $2^n$-dimensional space. As we did previously, we
index coordinates by subsets of $[n]$, so that a generic
vector is $(x_S)_{S\in[n]}$. In particular,
given a subset $T$ of $[n]$, the symbol
$e_T$ denotes the standard unit vector with a 1 in the $T$-th
coordinate and 0 elsewhere, i.e.,
$e_T=(\de_{S,T})_{S\in[n]}$. 

Now let $d,s$ be integers with $0\le d\le s \le n-d$ and let $A$ and
$B$ be
two disjoint $d$-element subsets $A$ and $B$ of $[n]$. Then we 
define the vector $v_{d,s}(A,B)$ by
\begin{equation}\label{e16}
v_{d,s}(A,B):=\underset {Y\subseteq [n]\backslash (A\cup B),\ 
\v{Y}=s-d}
{\sum _{X\subseteq A} ^{}}(-1)^{\v{X}}\, e_{X\cup X'\cup
Y},
\end{equation}
where $X'$ is the ``{\it complement of $X$ in $B$}" by which we mean that
if $X$ consists of the $i_1$-, $i_2$-, \dots -largest elements of $A$,
$i_1<i_2<\cdots$, then $X'$ consists of all elements of $B$ {\it
except for} the $i_1$-, $i_2$-, \dots -largest elements of $B$.
For example, let $n=7$. Then the vector $v_{2,3}(\{1,3\},\{2,5\})$ is
given by
\begin{multline}\label{e17}
 e_{\{2,4,5\}}+ e_{\{2,5,6\}}+
e_{\{2,5,7\}}- e_{\{1,4,5\}}- e_{\{1,5,6\}}- e_{\{1,5,7\}}\\
-
e_{\{2,3,4\}}- e_{\{2,3,6\}}- e_{\{2,3,7\}}+ e_{\{1,3,4\}}+
e_{\{1,3,6\}}+ e_{\{1,3,7\}}.
\end{multline}
(In this special case, the possible subsets $X$ of $A=\{1,3\}$ 
in the sum in (\ref{e16}) are $\emptyset$, 
$\{1\}$, $\{3\}$, $\{1,3\}$, with corresponding complements in
$B=\{2,5\}$ being $\{2,5\}$, $\{5\}$, $\{2\}$, $\emptyset$,
respectively,
and the possible sets $Y$ are $\{4\}$, $\{6\}$, $\{7\}$.) Observe
that all sets $X\cup X'\cup Y$ which occur as indices in (\ref{e16}) 
have the same cardinality $s$.

\begin{lemma}\label{l3} Let $d,s$ be integers with $0\le d\le s \le n-d$ and
let $A$ and $B$ be disjoint $d$-element subsets of $[n]$. Then
$v_{d,s}(A,B)$ as defined in {\em (\ref{e16})} is an eigenvector of the matrix
$\ze_n(u)$, the entries of which are given by {\em (\ref{e4})}, for
the eigenvalue $\la_d$, where $\la_d$ is given by {\em (\ref{e14})}.
\end{lemma}
We want to show that the multiplicity of $\la_d$ equals the
expression in (\ref{e15}). Of course, Lemma~\ref{l3} gives many more eigenvectors
for $\la_d$. Therefore, in order to describe a basis for the
corresponding eigenspace, we have to restrict the collection of
vectors in Lemma~\ref{l3}. 

We do this in the following way. Fix $d$, $0\le d\le \fl{n/2}$.
Let $P$ be a lattice path in the
plane integer lattice $\Z^2$, starting in $(0,0)$,
consisting of $n-d$ up-steps $(1,1)$ and
$d$ down-steps $(1,-1)$, which never goes below the $x$-axis.
Figure~1 displays an example with $n=7$ and $d=2$. Clearly, the end
point of $P$ is $(n,n-2d)$. We call a lattice path which starts in
$(0,0)$ and never goes below the $x$-axes a {\it ballot path}. (This
terminology is motivated by its relation to the (two-candidate) {\it ballot
problem\/}, see e.g\@. \cite[Ch.~1, Sec.~1]{MohaAE}. An alternative
term for ballot path which is often used is ``Dyck path'', see
e.g\@. \cite[p.~I-12]{vien}.) We will use the
abbreviation ``b.p.'' for ``ballot path'' in displayed formulas.
\vskip10pt
\vbox{
$$
\Gitter(8,5)(0,0)
\Koordinatenachsen(8,5)(0,0)
\Pfad(0,0),3433433\endPfad
\SPfad(7,0),222\endSPfad
\DickPunkt(0,0)
\DickPunkt(7,3)
\Label\l{\scriptstyle1}(1,1)
\Label\r{\scriptstyle2}(1,1)
\Label\l{\scriptstyle3}(3,1)
\Label\l{\scriptstyle4}(4,2)
\Label\r{\scriptstyle5}(4,2)
\Label\l{\scriptstyle6}(6,2)
\Label\l{\scriptstyle7}(7,3)
\hskip3.5cm
$$
\centerline{\small Ballot paths}
\vskip5pt
\centerline{\small Figure 1}
}
\vskip10pt

Given such a lattice path $P$, label the steps from $1$ to $n$, as
is indicated in Figure~1. Then
define $A_P$ to be set of all labels corresponding to the first $d$
up-steps of $P$ and $B_P$ to 
be set of all labels corresponding to the $d$ down-steps of $P$. In
the example of Figure~1 we have for the choice $d=2$ that 
$A_P=\{1,3\}$ and $B_P=\{2,5\}$.
Thus, to each $d$ and $s$, $0\le d\le s\le n-d$, and $P$ as above we
can associate the vector $v_{d,s}(A_P,B_P)$. In our running example
of Figure~1 the vector
$v_{2,3}(P)$ would hence be $v_{2,3}(\{1,3\},\{2,5\})$, the vector in
(\ref{e17}). To have a more concise form of
notation, we will write $v_{d,s}(P)$ for $v_{d,s}(A_P,B_P)$ from now on.

\begin{lemma}\label{l4}
The set of vectors
\begin{equation}\label{e18}
\{v_{d,s}(P):0\le d\le s\le n-d,\ P\text { a ballot path from $(0,0)$
to $(n,n-2d)$}\}
\end{equation}
is linearly independent.
\end{lemma}

The final Lemma tells us how many such vectors $v_{d,s}(P)$ there
are.
\begin{lemma} \label{l5}
The number of ballot paths from $(0,0)$ to $(n,n-2d)$
is $\frac {n-2d+1} {n+1}\binom {n+1}d$. The total number of all
vectors in the set\/ {\em(\ref{e18})} is $2^n$.
\end{lemma}

Now, let us for a moment assume that Lemmas~\ref{l3}--\ref{l5} are already proved.
Then, Theorem~\ref{t2} follows immediately, as it turns out.
\medskip

{\sc Proof of Theorem~\ref{t2}}. Consider the set of vectors in (\ref{e18}). By
Lemma~\ref{l3} we know that it consists of eigenvectors for the matrix
$\ze_n(u)$. In addition, Lemma~\ref{l4} tells us that this set of vectors is
linearly independent. Furthermore, by Lemma~\ref{l5} the number of vectors in
this set is exactly $2^n$, which is the dimension of the space where
all these vectors are contained. Therefore, they must form a basis of
the space. 

Lemma~\ref{l3} says more precisely that $v_{d,s}(P)$ is an eigenvector for
the eigenvalue $\la_d$. From what we already know, this implies that
for fixed $d$ the set
$$
\{v_{d,s}(P): d\le s\le n-d,\ P\text { a ballot path from $(0,0)$
to $(n,n-2d)$}\}$$
forms a basis for the eigenspace corresponding to $\la_d$. Therefore,
the dimension of the eigenspace corresponding to $\la_d$ equals
the number of possible numbers $s$ times the number of possible
lattice paths $P$. This is exactly
$$(n-2d+1)\frac {(n-2d+1)} {(n+1)}\binom {n+1}d,$$
the number of possible lattice paths $P$ being given by the first
statement of Lemma~\ref{l5}.
This expression equals exactly the expression (\ref{e15}).
Thus, Theorem~\ref{t2} is proved. \quad \quad \qed

\bigskip

Now we turn to the proofs of the Lemmas.
\medskip

{\sc Proof of Lemma~\ref{l3}}. Let $d,s$ and $A,B$ be fixed, satisfying the
restrictions in the statement of the Lemma.
We have to show that
$$\ze_n(u)\cdot v_{d,s}(A,B)=\la_d v_{d,s}(A,B).$$
Restricting our attention to the $I$-th component, we see 
from the
definition (\ref{e16}) of $v_{d,s}(A,B)$ that we need to establish 
\begin{equation}\label{e19}
\underset {Y\subseteq [n]\backslash (A\cup B),\ 
\v{Y}=s-d}
{\sum _{X\subseteq A} ^{}}\kern-1cm Z_{I,X\cup X'\cup Y}\,(-1)^{\v{X}}=
\begin{cases} \la_d (-1)^{\v{U}}&\text {if $I$ is of the form $U\cup U'\cup V$}\\
&\text {for some $U$ and $V$, $U\subseteq A$,}\\
&\text {$V\subseteq [n]\backslash (A\cup B)$,
$\v{V}=s-d$}\\
0&\text {otherwise.}\end{cases}
\end{equation}
We prove (\ref{e19}) by a case by case analysis. The first two cases cover
the case ``otherwise" in (\ref{e19}), the third case treats the first alternative
in (\ref{e19}).

\smallskip
{\it Case 1. The cardinality of $I$ is different from $s$}. As we
observed earlier, the cardinality of any set $X\cup X'\cup Y$ which
occurs as index at the left-hand side of (\ref{e19}) equals $s$. The
cardinality of $I$ however is different from $s$. 
As we observed in the Remark after Theorem~\ref{t1}, this implies that
any coefficient $Z_{I,X\cup X'\cup Y}$ on the left-hand side
vanishes. Thus, (\ref{e19}) is proved in this case.

\smallskip
{\it Case 2. The cardinality of $I$ equals $s$, but
$I$ does not have the form $U\cup U'\cup V$
for any $U$ and $V$, $U\subseteq A$,
$V\subseteq [n]\backslash (A\cup B)$, $\v{V}=s-d$}.
Now the sum on the left-hand side of (\ref{e19}) contains nonzero
contributions. We have to show that they cancel each other. We do
this by grouping summands in pairs, the sum of each pair being 0.

Consider a set $X\cup X'\cup Y$ which occurs as index at the
left-hand side of (\ref{e19}). Let $e$ be minimal such that 
\begin{enumerate}
\item[]either: the $e$-th largest element of $A$ and the
$e$-th largest element of $B$ are both in $I$,
\item[]or: the $e$-th largest element of $A$ and the
$e$-th largest element of $B$ are both not in $I$.
\end{enumerate}
That such an $e$ must exist is guaranteed by our assumptions about
$I$. Now consider $X$ and $X'$. If the $e$-th largest element of $A$
is contained in $X$ then the $e$-th largest element of $B$ is not
contained in $X'$, and vice versa. Define a new set $\bX$ by adding
to $X$ the $e$-th largest element of $A$ if it is not already
contained in $X$, respectively by removing it from $X$ if it is
contained in $X$. Then, it is easily checked that 
$$Z_{I,X\cup X'\cup Y}=Z_{I,\bX\cup \bX'\cup Y}.$$
On the other hand, we have $(-1)^{\v{X}}=-(-1)^{\v\bX}$ since the
cardinalities of $X$ and $\bX$ differ by $\pm1$. Both facts combined
give
$$Z_{I,X\cup X'\cup Y}\,(-1)^{\v{X}}+Z_{I,\bX\cup \bX'\cup Y}\,
(-1)^{\v\bX}=0.$$
Hence, we have found two summands on the left-hand side of
(\ref{e19}) which cancel each other. 

Summarizing,
this construction finds for any $X,Y$ sets $\bX,Y$ such that the
corresponding summands on the left-hand side of (\ref{e19}) cancel each
other. Moreover, this construction applied to $\bX,Y$ gives back
$X,Y$. Hence, what the construction does is exactly what we claimed,
namely it groups the summands into pairs which contribute 0 to the
whole sum. Therefore the sum is 0, which establishes (\ref{e19}) in this case
also.

\smallskip
{\it Case 3. $I$ has the form $U\cup U'\cup V$
for some $U$ and $V$, $U\subseteq A$,
$V\subseteq [n]\backslash (A\cup B)$, $\v{V}=s-d$}.
This assumption implies in particular that the cardinality of $I$ is
$s$. From the Remark after the statement of Theorem~\ref{t1} we know that in
our situation $Z_{I,X\cup X'\cup Y}$ depends only on the number of
common elements in $I$ and $X\cup X'\cup Y$. Thus, the left-hand side
in (\ref{e19}) reduces to
\begin{equation}\label{e20}
\sum _{j,k\ge0} ^{}N(j,k)\,(-1)^{\v{U}+j}\,k!
\frac {1} {2^n}\frac {\Ga\(\frac {5} {2}-u\)\, \Ga(2+n-s-u)\,
\Ga(2+s-u)} {\Ga\(\frac {5} {2}+\frac {n} {2}-u\)\, \Ga\(2+\frac {n}
{2}-u\)\Ga(2+k-u)},
\end{equation}
where $N(j,k)$ is the number of sets $X\cup X'\cup Y$,
for some $X$ and $Y$, $X\subseteq A$,
$Y\subseteq [n]\backslash (A\cup B)$, $\v{Y}=s-d$, which have $s-k$
elements in common with $I$, and which have $d-j$ elements in common
with $I\cap(A\cup B)=U\cup U'$. Clearly, we used expression (\ref{e4}) with
$n_{\in\in}=s-k$ and $n_{\notin\notin}=n-s-k$.

To determine $N(j,k)$, note first that there are $\binom dj$
possible sets $X\cup X'$ which intersect $U\cup U'$ 
in exactly $d-j$ elements. Next, let us assume that we already made a
choice for $X\cup X'$.
In order to determine the number of possible sets $Y$ such that
$X\cup X'\cup Y$ has $s-k$ elements in common with $I$, we have to
choose $(s-k)-(d-j)=s-d+j-k$ elements from
$V$, for which we have $\binom {s-d}{s-d+j-k}$ possibilities, and
we have to choose $s-d-(s-d+j-k)=k-j$ elements from $[n]\backslash
(I\cup A\cup B)$ to obtain a total number of $s$ elements, for which
we have $\binom {n-s-d}{k-j}$ possibilities. Hence,
\begin{equation}\label{e21}
N(j,k)=\binom dj\binom {s-d}{k-j}\binom {n-s-d}{k-j}.
\end{equation}

So it remains to evaluate the double sum (\ref{e20}), using the
expression (\ref{e21})
for $N(j,k)$. 

We start by writing the sum over $j$ in (\ref{e20}) in hypergeometric
notation,
\begin{multline}\notag 
(-1)^{\v{U}}\frac {1} {2^n}
{\frac{\Gamma({ \textstyle {\frac5 2} - u})  \,
     \Gamma({ \textstyle 2 + n - s - u})  \,\Gamma({ \textstyle 2 + s - u})  
}
{\Gamma({ \textstyle 2 - u})  \,
     \Gamma({ \textstyle 2 + {\frac n 2} - u})  \,
     \Gamma({ \textstyle {\frac5 2} + {\frac n 2} - u}) }}\\
\times
{     \sum_{k = 0}^{\infty}
        {\frac{  
            ({ \textstyle d - s}) _{k} \, ({ \textstyle d - n + s}) _{k} }
            {({ \textstyle 1}) _{k} \, ({ \textstyle 2 - u}) _{k} }}
{} _{3} F _{2} \!\left [ \begin{matrix} { -k, -k, -d}\\ { 1 - d - k + s, 1
             - d - k + n - s}\end{matrix} ; {\displaystyle 1}\right ]
}.
\end{multline}
To the $_3F_2$ series we apply a transformation formula of Thomae
(see e.g\@. \cite[(3.1.1)]{GaRaAA}),
\begin{equation}\label{e22}
{} _{3} F _{2} \!\left [ \begin{matrix} { a, b, -m}\\ { d, e}\end{matrix} ;
   {\displaystyle 1}\right ]  =
{\frac{ ({ \textstyle -b + e}) _{m} }
      {({ \textstyle e}) _{m} }}  
{} _{3} F _{2} \!\left [ \begin{matrix} { -m, b, -a + d}\\ { d, 1 + b - e -
       m}\end{matrix} ; {\displaystyle 1}\right ] 
\end{equation}
where $m$ is a nonnegative integer. We write the resulting $_3F_2$
again as a sum over $j$, then interchange sums over $k$ and $j$, and
write the (now) inner sum over $k$ in hypergeometric notation. Thus
we obtain
\begin{multline}\notag 
(-1)^{\v{U}}\frac {1} {2^n}
{\frac{\Gamma({ \textstyle {\frac5 2} - u})  \,
     \Gamma({ \textstyle 2 + n - s - u})  \,\Gamma({ \textstyle 2 + s - u})  
}
{    \Gamma({ \textstyle {\frac5 2} + {\frac n 2} - u}) \,
     \Gamma({ \textstyle 2 + {\frac n 2} - u})  \,
   \Gamma({ \textstyle 2 - u})}}\\
\times{ \sum_{j = 0}^{\infty}
        {\frac{  
            ({ \textstyle -d}) _{j} \, ({ \textstyle 1 - d + s}) _{j} } 
          {({ \textstyle 1}) _{j} \, ({ \textstyle 2 - u}) _{j} }}  }
{} _{2} F _{1} \!\left [ \begin{matrix} { j - n + s, d - s}\\ { 2 + j -
             u}\end{matrix} ; {\displaystyle 1}\right ].
\end{multline}
The $_2F_1$ series in this expression is terminating because $d-s$ is
a nonpositive integer. Hence, it can be summed by means of Gau\ss' sum
(\ref{e12}).
Writing the remaining sum over $j$ in hypergeometric notation, the
above expression becomes
$$(-1)^{\v{U}}\frac {1} {2^n}
{\frac{  
     \Gamma({ \textstyle {\frac5 2} - u})  \,
     \Gamma({ \textstyle 2 + n - d - u}) \, \Gamma({ \textstyle 2 + s - u}) }
     {      \Gamma({ \textstyle {\frac5 2} + {\frac n 2} - u})  \,
\Gamma({ \textstyle 2 + {\frac n 2} - u})  \,
     \Gamma({ \textstyle 2 + s -d - u}) }}
{} _{2} F _{1} \!\left [ \begin{matrix} { -d, 1 - d + s}\\ { 2 - d + s -
      u}\end{matrix} ; {\displaystyle 1}\right ].
$$
Again, 
the $_2F_1$ series is terminating and so
is summable by means of (\ref{e12}). Thus, we get
$$(-1)^{\v{U}}\frac {1} {2^n}\frac {\Ga\(\frac {5} {2}-u\)\, \Ga(2+n-d-u)\,
\Ga(1+d-u)} {\Ga\(\frac {5} {2}+\frac {n} {2}-u\)\, \Ga(2+\frac {n}
{2}-u)\,\Ga(1-u)},$$
which is exactly the expression (\ref{e14}) for $\la_d$ times $(-1)^{\v{U}}$. 
This proves (\ref{e19}) in this case.

\smallskip
The proof of Lemma~\ref{l3} is now complete.\quad \quad \qed
\medskip

\medskip
{\sc Proof of Lemma~\ref{l4}}.
We know from Lemma~\ref{l3} that $v_{d,s}(P)$ lies in the eigenspace for the
eigenvalue $\la_d$, with $\la_d$ being given in (\ref{e14}). The $\la_d$'s,
$d=0,1,\dots, \fl{n/2}$, are all distinct, so the corresponding
eigenspaces are linearly independent. Therefore it suffices to show
that for any {\it fixed\/} $d$ the set of vectors
$$\{v_{d,s}(P): d\le s\le n-d,\ P\text { a ballot path from $(0,0)$
to $(n,n-2d)$}\}$$
is linearly independent.

On the other hand, a vector $v_{d,s}(A,B)$ lies in the space spanned
by the standard unit vectors $e_T$ with $\v{T}=s$. Clearly, as $s$
varies, these spaces are linearly independent. Therefore, it suffices to
show that for any {\it fixed\/} $d$ {\it and\/} $s$ the set of vectors
$$\{v_{d,s}(P): P\text { a ballot path from $(0,0)$
to $(n,n-2d)$}\}$$
is linearly independent.

So, let us fix integers $d$ and $s$ with $0\le d\le s\le n-d$, and
let us suppose that there is some vanishing linear combination
\begin{equation}\label{e23}
\sum _{P\text{ b.p\@. from }(0,0)\text{ to }(n,n-2d)} ^{}
c_P\,v_{d,s}(P)=0.
\end{equation}
We have to establish that $c_P=0$ for all ballot paths $P$ from
$(0,0)$ to $(n,n-2d)$. 

We prove this fact by induction on the set of ballot paths from
$(0,0)$ to $(n,n-2d)$. In order to make this more precise, we need to
impose a certain order on the ballot paths. Given a ballot path $P$
from $(0,0)$ to $(n,n-2d)$, we define its {\it front portion\/} $F_P$
to be the portion of $P$ from the beginning up to and including $P$'s
$d$-th up-step. For example, choosing $d=2$, the front portion of the
ballot path in Figure~1 is the subpath from $(0,0)$ to $(3,1)$. Note
that  $F_P$ can be any ballot path starting in $(0,0)$ with
$d$ up-steps and less than $d$ down-steps. We order such front
portions lexicographically, in the sense that $F_1$ is before $F_2$ if
and only if $F_1$ and $F_2$ agree up to some point and then $F_1$
continues with an up-step while $F_2$ continues with a down-step.

Now, here is what we are going to prove: Fix any possible front
portion $F$. We shall show that $c_P=0$ for all $P$ with front portion
$F_P$ equal to $F$, {\it given that it is already known that
$c_{P'}=0$ for all $P'$ with a front portion $F_{P'}$ that is before
$F$.} Clearly, by induction, this would prove $c_P=0$ for {\it all\/}
ballot paths $P$ from $(0,0)$ to $(n,n-2d)$.

Let $F$ be a possible front portion, i.e., a ballot path starting in
$(0,0)$ with exactly $d$ up-steps and less than $d$ down-steps. As we
did earlier, label the steps of $F$ by $1,2,\dots$, and denote the set
of labels corresponding to the down-steps of $F$ by $B_F$. We write
$b$ for $\v{B_F}$, the number of all down-steps of $F$. Observe that
then the total number of steps of $F$ is $d+b$. 

Now, let $T$ be a fixed $(d-b)$-element subset of
$\{d+b+1,d+b+2,\dots,n\}$. Furthermore, let $S$ be a set of the form
$S=B_F\cup S_1\cup S_2$, where $S_1\subseteq T$ and $S_2\subseteq
\{d+b+1,d+b+2,\dots,n\}\backslash T$, and such that $\v{S}=s$. 

We consider the coefficient of $e_S$ in the left-hand side of (\ref{e23}). To
determine this coefficient, we have to determine the coefficient of
$e_S$ in $v_{d,s}(P)$, for all $P$. We may concentrate on those $P$
whose front portion $F_P$ is equal to or later than $F$, since our
induction hypothesis says that $c_P=0$ for all $P$ with $F_P$ before
$F$. So, let $P$ be a ballot path from $(0,0)$ to $(n,n-2d)$ with
front portion equal to or later than $F$. We claim that the
coefficient of $e_S$ in $v_{d,s}(P)$ is zero unless the set $B_P$ of
down-steps of $P$ is contained in $S$.

Let the coefficient of $e_S$ in $v_{d,s}(P)$ be nonzero.
To establish the claim, we first prove that the front portion $F_P$
of $P$ has to equal $F$. Suppose that this is not the case. Then the
front portion of $P$ runs in parallel with $F$ for some time, say for
the first $(m-1)$ steps, with some $m\le d+b$, and then $F$ continues
with an up-step and $F_P$ continues with a down-step (recall that
$F_P$ is equal to or later than $F$). By (\ref{e16}) we have
\begin{equation}\label{e24}
v_{d,s}(P):=\underset {Y\subseteq [n]\backslash (A_P\cup B_P),\ 
\v{Y}=s-d}
{\sum _{X\subseteq A_P} ^{}}(-1)^{\v{X}}\, e_{X\cup X'\cup
Y}.
\end{equation}
We are assuming that the coefficient of $e_S$ in $v_{d,s}(P)$ is
nonzero, therefore $S$ must be of the form $S=X\cup X'\cup Y$, with
$X,Y$ as described in (\ref{e24}). We are considering the case that the
$m$-th step of $F_P$ is a down-step, whence $m\in B_P$, while the
$m$-th step of $F$ is an up-step, whence $m\notin B_F$. By definition
of $S$, we have $S\cap \{1,2\dots,d+b\}=B_F$, whence $m\notin S$. 

Summarizing so far, we have $m\in B_P$, $m\notin S$, for some $m\le
d+b$, and $S=X\cup X'\cup Y$, for some $X,Y$ as described in (\ref{e24}). In
particular we have $m\notin X'$. Now recall that $X'$ is the
``complement of $X$ in $B_P$''. This says in particular that, if $m$ is
the $i$-th largest element in $B_P$, then the $i$-th largest element
of $A_P$, $a$ say, is an element of $X$, and so of $S$. By
construction of $A_P$ and $B_P$, $a$ is smaller than $m$, so in
particular $a<d+b$. As we already observed, there holds $S\cap
\{1,2,\dots,d+b\}=B_F$, so we have $a\in B_F$, i.e., the $a$-th step
of $F$ is a down-step. On the other hand, we assumed that $P$ and $F$
run in parallel for the first $(m-1)$ steps. Since $a\in A_P$, the set
of up-steps of $P$, the $a$-th step of $P$ is an up-step. We have
$a\le m-1$, therefore the $a$-th step of $F$ must be an up-step
also. This is absurd. Therefore, given that the
coefficient of $e_S$ in $v_{d,s}(P)$ is nonzero, the front portion
$F_P$ of $P$ has to equal $F$.

Now, let $P$ be a ballot path from $(0,0)$ to $(n,n-2d)$ with
front portion equal to $F$, and suppose that $S$ has the form $S=X\cup
X'\cup Y$, for some $X,Y$ as described in (\ref{e24}). By definition of the
front portion, the set $A_P$ of up-steps of $P$ has the property
$A_P\cap \{1,2,\dots,d+b\}=\{1,2,\dots,d+b\}\backslash B_F$. Since
$\v{B_F}=b$, these are the labels of exactly $d$ up-steps. Since the
cardinality of $A_P$ is exactly $d$ by definition, we must have
$A_P=\{1,2,\dots,d+b\}\backslash B_F$. Because of $S\cap
\{1,2,\dots,d+b\}=B_F$, which we already used a number of times, $A_P$
and $S$ are disjoint, which in particular implies that $A_P$ and $X$
are disjoint. However, $X$ is a subset of $A_P$ by definition, so $X$
must be empty. This in turn implies that $X'=B_P$. This says nothing
else but that the set $B_P$ of down-steps of $P$ equals $X'$ and so is
contained in $S$. This establishes our claim. 

In fact, we proved more. We saw that $S$ has the form $S=X\cup X'\cup
Y$, with $X=\emptyset$. This implies that the coefficient of $e_S$ in
$v_{d,s}(P)$, as given by (\ref{e24}), is actually $+1$. Comparison of
coefficients of $e_S$ in (\ref{e23}) then gives
\begin{equation}\label{e25}
\underset {F_P=F,\ B_P\subseteq S}
{\sum _{P\text{ b.p\@. from }(0,0)\text{ to }(n,n-2d)} ^{}}
c_P=0,
\end{equation}
for any $S=B_F\cup S_1\cup S_2$, where $S_1\subseteq T$ and $S_2\subseteq
\{d+b+1,d+b+2,\dots,n\}\backslash T$, and such that $\v{S}=s$.

Now, we sum both sides of (\ref{e25}) over all such sets $S$, keeping the
cardinality of $S_1$ and $S_2$ fixed, say
$\v{S_1}=d-b-j$, enforcing $\v{S_2}=s-d+j$, for a fixed $j$, $0\le
j\le d-b$. For a fixed ballot path $P$ from $(0,0)$ to $(n,n-2d)$,
with front portion $F$, with $d-b-k$ down-steps in $T$, and hence with
$k$ down-steps in $\{d+b+1,d+b+2,\dots,n\}\backslash T$, there are
$\binom k{k-j}$ such sets $S_1\subseteq T$ containing all the
$d-b-k$ down-steps of $P$ in $T$, and there are $\binom
{n-(d+b)-(d-b)-k} {s-d+j-k}$ such sets $S_2\subseteq
\{d+b+1,d+b+2,\dots,n\}\backslash T$ containing all the $k$ down-steps
of $P$ in $\{d+b+1,d+b+2,\dots,n\}\backslash T$. Therefore,
summing up (\ref{e25}) gives
\begin{equation}\label{e26}
\sum _{k\ge0} ^{}\binom kj\binom {n-2d-k}{n-d-s-j} 
\bigg(\kern-5pt\underset {\v{B_P\cap(\{d+b+1,d+b+2,\dots,n\}\backslash T)}=k}
{\underset {F_P=F,\ \v{B_P\cap T}=d-b-k}
{\sum _{P\text{ b.p\@. from }(0,0)\text{ to }(n,n-2d)} ^{}}}\kern-.3cm
c_P\bigg)=0,\quad 
j=0,1,\dots,d-b.
\end{equation}
Denoting the inner sum in (\ref{e26}) by $C(k)$, we see that
(\ref{e26}) represents
a non-degenerate triangular system of linear equations for $C(0),
C(1),\dots, C(d-b)$. Therefore, all the quantities $C(0),
C(1),\dots, C(d-b)$ have to equal 0. In particular, we have
$C(0)=0$. Now, $C(0)$ consists of just a single term $c_P$, with $P$
being the ballot path from $(0,0)$ to $(n,n-2d)$, with front portion
$F$, and the labels of the $d-b$ down-steps besides those of $F$
being exactly the elements of $T$. Therefore, we have $c_P=0$ for this ballot
path. The set $T$ was an arbitrary $(d-b)$-subset of
$\{d+b+1,d+b+2,\dots,n\}$. Thus, we have proved $c_P=0$ for any ballot
path $P$ from $(0,0)$ to $(n,n-2d)$ with front portion $F$. This
completes our induction proof. \quad \quad \qed

\medskip

\medskip
{\sc Proof of Lemma~\ref{l5}}. That the number of ballot paths from $(0,0)$
to $(n,n-2d)$ equals $\frac {n-2d+1} {n+1}\binom {n+1}d$ is a
classical combinatorial result (see e.g\@. \cite[Theorem~1
with $t=1$]{MohaAE}). From
this it follows that the total number of vectors in the set (\ref{e18}) is
\begin{equation}\label{e27}
\sum _{d=0} ^{\fl{n/2}}(n-2d+1)\frac {(n-2d+1)} {(n+1)}\binom
{n+1}d.
\end{equation}
To evaluate this sum, note that the summand is invariant under the
substitution $d\to n-2d+1$. Therefore, extending the range of
summation in (\ref{e27}) to $d=0,1,\dots,n+1$ and dividing the result by $2$
gives the same value. So, the cardinality of the set (\ref{e18}) is also
given by
$$\frac {1} {2}\sum _{d=0} ^{n+1}\frac {(n-2d+1)^2} {(n+1)}\binom
{n+1}d.$$
Using the simple identity
$$\frac {(n-2d+1)^2} {(n+1)}\binom {n+1}d=
(n+1)\binom {n+1}d-4n\binom n{d-1}+4n \binom {n-1}{d-2},$$
the last sum can be decomposed into 
$$
\frac {n+1} {2}\sum _{d=0} ^{n+1}\binom {n+1}d-
2n\sum _{d=1} ^{n+1}\binom 
n{d-1}+2n\sum _{d=2} ^{n+1} \binom {n-1}{d-2}.
$$
Each of these sums can be evaluated by the binomial theorem, and thus
the expression reduces to $2^n$. This completes the proof of the Lemma.\quad \quad
\qed

\medskip

In fact, Theorem~\ref{t2} can be generalized to a wider class of matrices. 

\begin{theorem}\label{t6}
Let $\tilde \ze_n(u)=(\tilde Z_{IJ})_{I,J\in [n]}$ be the $2^n\times
2^n$ matrix defined by
$$\tilde Z_{IJ}:=\de_{n_{\notin\in},n_{\in\notin}} \frac {\big(\tfrac
{n-n_{\in\in}- n_{\notin\notin}} 2\big)! } {\Ga\(2+
\tfrac
{n-n_{\in\in}- n_{\notin\notin}} 2-u\)}\cdot
f(n_{\in\in}-n_{\notin\notin}),$$
where $n_{\in\in}$, etc., have the same meaning as earlier, and where
$f(x)$ is a function of $x$ which is symmetric, i.e., $f(x)=f(-x)$.
Then, the eigenvalues of 
$\tilde\ze_n(u)$ are
\begin{equation}\label{e28}
\la_{d,s}=f(n-2s)\frac { \Ga(2+n-d-u)\,
\Ga(1+d-u)} {\Ga(2+n-s-u) \,\Ga(2+s-u)\,
\Ga(1-u)}, \quad 0\le d\le s\le n-d,
\end{equation}
with respective multiplicities
\begin{equation}\label{e29}
\frac {n-2d+1} {n+1}\binom {n+1}d,
\end{equation}
independent of $s$.

\end{theorem}
\medskip

{\sc Proof}.
The above proof of Theorem~\ref{t2} has to be adjusted only
insignificantly to yield a proof of Theorem~\ref{t6}. In particular, the
vector $v_{d,s}(A,B)$ as defined in (\ref{e16}) is an eigenvector for
$\la_{d,s}$, for any two disjoint $d$-element subsets $A$ and $B$ 
of $[n]$, and
the set (\ref{e18}) is a basis of eigenvectors for $\tilde\ze_n(u)$.
\quad \quad \qed
\medskip

\subsection{The relative entropies of $
 \raise6pt\hbox{$^n$}\kern-9pt \otimes\rho$ 
with
respect to the Bayesian density matrices $\ze_n(u)$}\label{s2.3}

We now apply the preceding results to compute the relative entropy of
$\overset n\otimes \rho$ with respect to $\ze_n(u)$. Utilizing the
definition (\ref{eq:5}) of relative entropy and employing the property
\cite{oh,wehrl} that $S(\overset n\otimes\rho)=nS(\rho)$,
it is given by
\begin{equation}\label{e31}
-n\,S(\rho)-\Tr\(\overset
n\otimes \rho\cdot \log\ze_n(u)\).
\end{equation}
The term $S(\rho)$ has been given in (\ref{eq:7}). Concerning the
second term in (\ref{e31}), we have the following theorem.
\begin{theorem}\label{t7}
Let $\ze_n(u)=(Z_{IJ})_{I,J\in[n]}$ be the matrix with
entries $Z_{IJ}$ given in {\em(\ref{e4})}. Then, we have
\begin{multline}\label{e30}
\Tr\(\overset n\otimes \rho\cdot\log \ze_n(u)\)\\
=\sum _{d=0}
^{\fl{n/2}} \frac {n-2d+1} {n+1}\binom {n+1}d\frac {1} {2^{n+1}r}
\big((1+r)^{n+1-d}(1-r)^d-(1+r)^d(1-r)^{n+1-d}\big)\log\la_d,
\end{multline} 
with $\la_d$ as given in {\em(\ref{e14})}.
\end{theorem}
Before we move on to the proof, we note that Theorem~\ref{t7} gives us 
the following expression for the relative entropy of 
$\overset n\otimes \rho$ with respect to $\ze_n(u)$

\begin{corollary}\label{c8}
The relative entropy of 
$\overset n\otimes \rho$ with respect to $\ze_n(u)$ equals
\begin{multline}\label{e32}   
\frac {n} {2}(1-r)\log((1-r)/2)+\frac {n}
{2}(1+r)\log((1+r)/2)\\
-\sum _{d=0} ^{\fl{n/2}}\frac {(n-2d+1)} {(n+1)}\binom {n+1}d\hskip5cm\\
\cdot\frac
{1} {2^{n+1}r}\((1+r)^{n-d+1}(1-r)^d-(1+r)^d(1-r)^{n-d+1}\)\log\la_d,
\end{multline}
with $\la_d$ as given in {\em(\ref{e14})}. 
\end{corollary}

\medskip
{\sc Proof of Theorem~\ref{t7}}.
One way of determining the trace of a linear
operator $L$ is 
to choose a basis of the vector space, $\{v_I:I\in [n]\}$ say,
write the action of $L$ on the basis elements in the form
$$Lv_I=c_Iv_I+\text{linear combination of $v_J$'s, $J\ne I$},$$
and then form the sum $\sum _{I} ^{}c_I$ of the ``diagonal'' 
coefficients, which
gives exactly the trace of $L$.

Clearly, we choose as a basis our set (\ref{e18}) of eigenvectors for
$\ze_n(u)$. To determine the action of $\overset n\otimes\rho
\cdot\log\ze_n(u)$ we need only to find the action of $\overset
n\otimes\rho$ on the vectors in the set (\ref{e18}). We claim that this
action can be described as
\begin{multline}\label{e33}
\big(\overset n\otimes\rho\big)\cdot v_{d,s}(P)\\
=\frac {1} {2^n}\bigg(\sum _{k\ge j\ge0} ^{}
(-1)^j\binom dj\binom {s-d}{k-j}\binom {n-s-d}{k-j}(1+z)^{s-k}
(x^2+y^2)^k(1-z)^{n-s-k}\bigg)\\
\cdot v_{d,s}(P)\ {}+{}\ \text{linear combination of eigenvectors}\\
\text{$v_{d',s'}(P')$ with $s'\ne s$},
\end{multline}
for any basis vector $v_{d,s}(P)$ in (\ref{e18}).

To see this, consider the $I$-th component of $\big(\overset n\otimes
\rho\big)\cdot v_{d,s}(P)$, i.e., the coefficient of $e_I$ in
$\big(\overset n\otimes
\rho\big)\cdot v_{d,s}(P)$, $I\in[n]$. By the definition (\ref{e16}) of
$v_{d,s}(P)$ it equals
\begin{equation}\label{e34}
\underset {Y\subseteq [n]\backslash (A_P\cup B_P),\ 
\v{Y}=s-d}
{\sum _{X\subseteq A_P} ^{}}\kern-1cm R_{I,X\cup X'\cup
Y}\,(-1)^{\v{X}},
\end{equation}
where $R_{IJ}$ denotes the $(I,J)$-entry of $\overset
n\otimes\rho$. (Recall that $R_{IJ}$ is given explicitly in (\ref{e2}).)
Now, it should be observed that we did a similar calculation already,
namely in the proof of Lemma~\ref{l3}. In fact, the expression (\ref{e34}) is almost
identical with the left-hand side of (\ref{e19}). The essential difference is
that $Z_{IJ}$ is replaced by $R_{IJ}$ for all $J$ (the nonessential
difference is that $A,B$ are replaced by $A_P,B_P$, respectively).  
Therefore, we can partially rely upon
what was done in the proof of Lemma~\ref{l3}. 

We distinguish between the same cases as in the proof of Lemma~\ref{l3}.

\smallskip
{\it Case 1. The cardinality of $I$ is different from $s$}. We do not
have to worry about this case, since $e_I$ then lies in the span of
vectors $v_{d',s'}(P')$ with $s'\ne s$, which is taken care of in
(\ref{e33}). 

\smallskip
{\it Case 2. The cardinality of $I$ equals $s$, but
$I$ does not have the form $U\cup U'\cup V$
for any $U$ and $V$, $U\subseteq A_P$,
$V\subseteq [n]\backslash (A_P\cup B_P)$, $\v{V}=s-d$}.
Essentially the same arguments as those in Case~2 in the proof of
Lemma~\ref{l3} show that the term (\ref{e34}) vanishes for this choice of $I$. Of
course, one has to use the explicit expression (\ref{e2}) for $R_{IJ}$.

\smallskip
{\it Case 3. $I$ has the form $U\cup U'\cup V$
for some $U$ and $V$, $U\subseteq A_P$,
$V\subseteq [n]\backslash (A_P\cup B_P)$, $\v{V}=s-d$}.
In Case~3 in the proof of Lemma~\ref{l3} we observed that there are
 $N(j,k)$ sets $X\cup X'\cup Y$,
for some $X$ and $Y$, $X\subseteq A_P$,
$Y\subseteq [n]\backslash (A_P\cup B_P)$, $\v{Y}=s-d$, which have $s-k$
elements in common with $I$, and which have $d-j$ elements in common
with $I\cap(A_P\cup B_P)=U\cup U'$, where $N(j,k)$ is given by
(\ref{e21}). Then, using the explicit expression (\ref{e2}) for $R_{IJ}$,
it is straightforward to see that the expression (\ref{e34}) equals
$$\frac {1} {2^n}\sum _{k\ge j\ge0} ^{}
(-1)^{\v U+j}\binom dj\binom {s-d}{k-j}\binom {n-s-d}{k-j}(1+z)^{s-k}
(x^2+y^2)^k(1-z)^{n-s-k}$$
in this case. This establishes (\ref{e33}).

\smallskip
Now we are in the position to write down an expression for the trace
of $\overset n\otimes\rho\cdot\log\ze_n(u)$. By Theorem~\ref{t2} and by (\ref{e33})
we have
\begin{multline}\label{e35}
\(\overset n\otimes\rho\cdot \log\ze_n(u)\)
\cdot v_{d,s}(P)\\
=\frac {1} {2^n}\bigg(\sum _{k\ge j\ge0} ^{}
(-1)^j\binom dj\binom {s-d}{k-j}\binom {n-s-d}{k-j}(1+z)^{s-k}
(x^2+y^2)^k(1-z)^{n-s-k}\bigg)\\
\cdot \log\la_d\cdot v_{d,s}(P)+\text{linear combination of eigenvectors}\\
\text{$v_{d',s'}(P')$ with $s'\ne s$}.
\end{multline}
 From what was said at the beginning of this proof, in order
to obtain the trace of $\overset n\otimes\rho\cdot\log\ze_n(u)$, we have to
form the sum of all the ``diagonal'' coefficients in (\ref{e35}).
Using the first statement of Lemma~\ref{l5} and replacing $x^2+y^2$ by
$r^2-z^2$, we see that it is
\begin{multline}\label{e36}
\sum _{d=0} ^{\fl{n/2}}\log\la_d \,\frac {(n-2d+1)} {(n+1)}\binom
{n+1}d\frac {1} {2^n}
\sum _{s=d} ^{n-d}\sum _{k\ge j\ge0} ^{}
(-1)^j\binom dj\binom {s-d}{k-j}\binom {n-s-d}{k-j}\\
\cdot (1+z)^{s-k}
(r^2-z^2)^k(1-z)^{n-s-k}.
\end{multline}
In order to see that this expression equals (\ref{e30}), we have to prove
\begin{multline}\label{e37} 
\sum _{s=d} ^{n-d}\sum _{j=0} ^{d}\sum _{k=j} ^{s}
(-1)^j\binom dj\binom {s-d}{k-j}\binom {n-s-d}{k-j} (1+z)^{s-k}
(r^2-z^2)^k(1-z)^{n-s-k}\\
=\frac {1} {2r}
\big((1+r)^{n+1-d}(1-r)^d-(1+r)^d(1-r)^{n+1-d}\big).
\end{multline}

We start with the left-hand side of (\ref{e37}) and write the inner sum in
hypergeometric notation, thus obtaining
$$\sum_{s = d}^{ n-d}\sum_{j = 0}^{d}
     {{\left( 1 - z \right) }^{ n - s-j}} 
         {{\left( 1 + z \right) }^{s-j}} 
         {{\left( {r^2} - {z^2} \right) }^j} 
        {\frac{     ({ \textstyle -d}) _{j} } {({ \textstyle 1}) _{j} }}
     {} _{2} F _{1} \!\left [ \begin{matrix} {  d - n + s,d-s}\\ { 1}\end{matrix}
        ; \frac {r^2-z^2} {1-z^2}\right].
$$
To the $_2F_1$ series we apply the transformation formula
(\cite[(1.8.10),
terminating form]{SlatAC}
$$
{} _{2} F _{1} \!\left [ \begin{matrix} { a, -m}\\ { c}\end{matrix} ; {\displaystyle
   z}\right ]  = 
{\frac{    ({ \textstyle c-a}) _{m} } {({ \textstyle c}) _{m} }}
{} _{2} F _{1} \!\left [ \begin{matrix} { -m, a}\\ { 1 + a - c -
       m}\end{matrix} ; {\displaystyle 1 - z}\right ] , 
$$
where $m$ is a nonnegative integer. We write the resulting $_2F_1$
series again as a sum over $k$. In the resulting expression
we exchange sums so that the sum over $j$ becomes the innermost sum.
Thus, we obtain
\begin{multline}\notag
 \sum_{s = d}^{n-d}\sum_{k = 0}^{s-d}
     {{\left( 1 - r^2 \right) }^k}
         {{\left( 1 - z \right) }^{n - s-k}} 
         {{\left( 1 + z \right) }^{s-k}} \\
\cdot       {\frac{  ({ \textstyle d - s}) _{k}\,  
         ({ \textstyle n - d - s+1}) _{s-d}\,  
         ({ \textstyle d - n + s}) _{k} } 
       { ({ \textstyle 1}) _{k}\,  
         ({ \textstyle 1}) _{s-d}  \,({ \textstyle 2 d - n}) _{k} }}
\sum _{j=0} ^{d}\binom dj \(\frac {z^2-r^2} {1-z^2}\)^j.
\end{multline}
Clearly, the innermost sum can be evaluated by the binomial theorem.
Then, we interchange sums over $s$ and $k$. The expression
that results is
\begin{multline}\notag \sum_{k = 0}^{\fl{n/2}-d }
{{\left( 1 - r^2 \right) }^{d + k}} 
       {{\left( 1 - z \right) }^{n-2d - 2k }} 
{\frac{       ({ \textstyle 2 d + k - n}) _{k} } {({ \textstyle 1}) _{k} }}\\
\cdot \sum _{s=0} ^{n-2d-2k}\binom {n-2d-2k}s \(\frac {1+z} {1-z}\)^s.
\end{multline}
Again, we can apply the binomial theorem.
Thus,  we reduce our expression on the left-hand
side of (\ref{e37}) to
$${2^{n-2 d}} {{\left( 1 - r^2 \right) }^d} 
 \sum _{k=0} ^{\fl{n/2}-d}\frac {\(d-\frac {n} {2}\)_k\, \(d-\frac
{n} {2}+\frac {1} {2}\)_k} {(2d-n)_k\, k!}(1-r^2)^k.$$
Now, we replace $(1-r^2)^k$ by its binomial expansion $\sum _{l=0}^k(-1)^l
\binom kl r^{2l}$, interchange sums over $k$ and $l$, and write
the (now) inner sum over $k$ in hypergeometric notation. This gives
\begin{multline}\notag {2^{n-2 d }} {{\left( 1 - r^2 \right) }^d} 
  \bigg( \sum_{l = 0}^{\fl{n/2}-d}{{\left( -1 \right) }^l} {r^{2 l}} 
{\frac{
    ({ \textstyle d - {\frac n 2}}) _{l}\,  
         ({ \textstyle {\frac1 2} + d - {\frac n 2}}) _{l} } 
       {({ \textstyle 1}) _{l}\,  ({ \textstyle 2 d - n}) _{l} }}  \\
     \cdot {} _{2} F _{1} \!\left [ \begin{matrix} { d + l - {\frac n
  2}, {\frac 1 2} + d
          + l - {\frac n 2}}\\ { 2 d + l - n}\end{matrix} ; {\displaystyle
          1}\right ]\bigg) .
\end{multline}
Finally, this $_2F_1$ series can be summed by means of Gau\ss'
summation (\ref{e12}). Simplifying, we have
$$ {{\left( 1 - r ^2\right) }^d}
  \sum _{l=0} ^{\fl{n/2}-d}\binom {n-2d+1}{2l+1}r^{2l},$$
which is easily seen to equal the right-hand side in (\ref{e37}). This
completes the proof of the Theorem. \quad \quad \qed

\medskip

\subsection{Asymptotics of the relative entropy
 of $\raise6pt\hbox{$^n$}\kern-9pt \otimes\rho$ with respect to 
$\ze_n(u)$}\label{s2.4}
In the preceding subsection, we obtained in Corollary~\ref{c8}
the general formula (\ref{e32}) for the relative entropy of $\overset
n\otimes\rho$ with respect to the Bayesian density matrix
$\ze_n(u)$. We, now, proceed to find its asymptotics for $n\to\infty$.
We prove the following theorem.

\begin{theorem}\label{t9}
The asymptotics of the relative entropy of 
$\overset n\otimes \rho$ with respect to $\ze_n(u)$ for a fixed $r$
with $0\le r<1$ is given by
\begin{multline}\label{a3} 
\frac {3} {2}\log n  -\frac {1} {2}- \frac {3} {2} \log 2  
-(1-u)\log(1-r^2)+\frac {1} {2r}\log\(\frac {1-r} {1+r}\)\\
+ \log\Ga(1 - u)  -\log\Ga(5/2 - u)
 +O\(\frac {1} {n}\).
\end{multline}
In the case $r=0$, this means that the asymptotics is given by the
expression {\em(\ref{a3})} in the limit $r\downarrow0$, i.e., by
\begin{equation}\label{a4} 
\frac {3} {2}\log n  -\frac {3} {2}- \frac {3} {2} \log 2  
+ \log\Ga(1 - u) -\log\Ga(5/2 - u) 
 +O\(\frac {1} {n}\).
\end{equation}
For any fixed $\ep>0$, the $O(.)$ term in {\em(\ref{a3})} 
is uniform in $u$ and $r$ as long as $0\le r\le 1-\ep$.

For $r=1$ the asymptotics is given by
\begin{equation}\label{a5}
 (2-u)\log n+(2u-3)\log 2
+ \frac {1} {2} \log\pi - \log\Ga(5/2 - u) + O\(\frac {1} {n}\).
\end{equation}
Also here, the $O(.)$ term is uniform in $u$.
\end{theorem}
\begin{remark}\em 
It is instructive to observe that, although a comparison of (\ref{a3})
and (\ref{a5}) seems to suggest that the
asymptotics of the relative entropy of 
$\overset n\otimes \rho$ with respect to $\ze_n(u)$ behaves completely 
differently for $0\le r<1$ and $r=1$, the two cases are really quite
compatible. In fact, letting $r$ tend to $1$ in (\ref{a3}) shows that
(ignoring the error term) the asymptotic expression approaches
$+\infty$ for $u<1/2$, $-\infty$ for $u>1/2$, and it approaches
$ \frac {3} {2}\log n-\frac {1} {2}-\frac {5} {2}\log 2
+ \frac {1} {2} \log\pi $ for $u=1/2$. This indicates that, for
$r=1$, the order of magnitude of the relative entropy of 
$\overset n\otimes \rho$ with respect to $\ze_n(u)$ should be larger
than $\frac {3} {2}\log n$ if $u<1/2$, smaller than $\frac {3}
{2}\log n$ if $u>1/2$, and exactly $\frac {3} {2}\log n$ if $u=1/2$.
How much larger or smaller is precisely what formula (\ref{a5})
tells us: the order of magnitude is $(2-u)\log n$, and in the case $u=1/2$
the asymptotics is, in fact, $ \frac {3} {2}\log n-2\log 2
+ \frac {1} {2} \log\pi $.
\end{remark}
\medskip

The proof of Theorem~\ref{t9} relies on several auxiliary
summations and estimations. These are
stated and proved separately in Lemma~\ref{la1} and \ref{la2}.

{\sc Proof of Theorem~\ref{t9}}. 
We start with the case $0<r<1$.
We concentrate first on the sum in (\ref{e32}).
Because of $\la_{n+1-d}=\la_d$, 
we have
\begin{multline}\notag
\frac {1} {2^{n+1}r}
\sum _{d=0} ^{\fl{n/2}}\frac {(n-2d+1)} {(n+1)}\binom {n+1}d\hskip5cm\\
\cdot\((1+r)^{n-d+1}(1-r)^d-(1+r)^d(1-r)^{n-d+1}\)\log\la_d\\
=\frac {1} {2^{n+1}r}\sum _{d=0} ^{n+1}\frac {(n-2d+1)} {(n+1)}\binom {n+1}d
(1+r)^{n-d+1}(1-r)^d\log\la_d.
\end{multline}
We expand the logarithm according to the addition rule to obtain 
\begin{align}\notag
 \frac {1} {2^{n+1}r}&
\sum _{d=0} ^{\fl{n/2}}\frac {(n-2d+1)} {(n+1)}\binom {n+1}d\\
\notag
&\hskip3cm\cdot\((1+r)^{n-d+1}(1-r)^d-(1+r)^d(1-r)^{n-d+1}\)\log\la_d\\
\notag
=&\frac {1} {2^{n+1}r}\sum _{d=0} ^{n+1}\frac {(n-2d+1)}
{(n+1)}\binom {n+1}d (1+r)^{n-d+1}(1-r)^d\\
\notag
&\hskip3cm\cdot\log\(\frac {1} {2^n}
\frac {\Ga(5/2-u)}
{\Ga(5/2+n/2-u)\,\Ga(2+n/2-u)\,\Ga(1-u)}\)\\
\notag
&\hskip.2cm +\frac {1} {2^{n+1}r}\sum _{d=0} ^{n+1}
{\frac {(n-2d+1)} {(n+1)}\binom
{n+1}d}(1+r)^{n-d+1}(1-r)^d
\log\Ga(1+d-u)\\
\label{a6}
&\hskip.2cm -\frac {1} {2^{n+1}r}\sum _{d=0} ^{n+1}
{\frac {(n-2d+1)} {(n+1)}\binom
{n+1}d}(1-r)^{n-d+1}(1+r)^d
\log\Ga(1+d-u).
\end{align}
The first sum on the right-hand side of (\ref{a6}) can be evaluated by means
of (\ref{a8}).
Besides, by Stirling's formula we have
$$\log\Ga(z)=\(z-\frac {1} {2}\)\log(z)-z+\frac {1} {2}\log2
+\frac {1} {2}\log\pi+O\(\frac {1} {z}\).$$
Thus, we get
\begin{multline}\label{a7}
\frac {1} {2^{n+1}r}
\sum _{d=0} ^{\fl{n/2}}\frac {(n-2d+1)} {(n+1)}\binom {n+1}d\\
\hskip3cm\cdot\((1+r)^{n-d+1}(1-r)^d-(1+r)^d(1-r)^{n-d+1}\)\log\la_d\\
=-n\log 2-
\log\Ga(5/2+n/2-u)-\log \Ga(2+n/2-u)+\log\Ga(5/2-u)\hskip30pt\\
\hskip.2cm -\log \Ga(1-u)
+\frac {1} {2^{n+1}r}\sum _{d=0} ^{n+1}
{\frac {(n-2d+1)} {(n+1)}\binom
{n+1}d}(1+r)^{n-d+1}(1-r)^d\hskip30pt\\
\hskip10pt\cdot\((1/2-u+d)\log(1+d-u)-(1-u+d)+\frac {1} {2}\log2
+\frac {1} {2}\log\pi+O\(\frac {1} {d+1}\)\)\\
-\frac {1} {2^{n+1}r}\sum _{d=0} ^{n+1}
{\frac {(n-2d+1)} {(n+1)}\binom
{n+1}d}(1-r)^{n-d+1}(1+r)^d\hskip3cm\\
\hskip10pt\cdot\((1/2-u+d)\log(1+d-u)-(1-u+d)+\frac {1} {2}\log2
+\frac {1} {2}\log\pi+O\(\frac {1} {d+1}\)\).\\
\end{multline}
Now, the sums are split into several sums by additivity. Those which
arise from the first sum in (\ref{a7}) can be evaluated using  
(\ref{a8}),
(\ref{a9}), (\ref{a10}), or approximated using (\ref{a15}). Those
which arise from the second sum can be evaluated by the same
identities and approximations, only with $r$ replaced by its negative.
Thus, we obtain
\begin{multline}\notag
\frac {1} {2^{n+1}r}
\sum _{d=0} ^{\fl{n/2}}\frac {(n-2d+1)} {(n+1)}\binom {n+1}d\hskip5cm\\
\cdot\((1+r)^{n-d+1}(1-r)^d-(1+r)^d(1-r)^{n-d+1}\)\log\la_d\\
=\frac {n} {2}(1-r)\log((1-r)/2)+
\frac {n} {2}(1+r)\log((1+r)/2)
-\frac {3} {2}\log n  + \frac {3} {2} \log 2 + \frac {1} {2}\\ 
+(1-u)\log(1-r)+(1-u)\log(1+r)+\frac {1} {2r}\log\(\frac {1+r} {1-r}\)\\
+ \log\Ga(5/2 - u)- \log\Ga(1 - u)  
 +O\(\frac {1} {n}\).
\end{multline}
Finally, use of this in (\ref{e32}) gives the claimed asymptotics
(\ref{a3}) for the
relative entropy of 
$\overset n\otimes \rho$ with respect to $\ze_n(u)$.

A closer analysis of the error terms shows that they can, in fact, be
bounded uniformly in $u$ and $r$, $0<r\le 1-\ep$, for any fixed
positive $\ep$.

\medskip
Now we turn to the two exceptional cases $r=0$ and $r=1$. 

In the case $r=1$, by (\ref{e32}) the relative entropy of 
$\overset n\otimes \rho$ with respect to $\ze_n(u)$ equals
$$
\frac {n} {2}(1-r)\log((1-r)/2)+\frac {n}
{2}(1+r)\log((1+r)/2)-\log\la_0,
$$
$\la_0$ being given by (\ref{e14}). A straightforward application of
Stirling's formula then leads to (\ref{a5}).

In the case $r=0$, the relative entropy (\ref{e32}) of 
$\overset n\otimes \rho$ with respect to $\ze_n(u)$ reduces to
\begin{multline}\notag
\frac {n} {2}(1-r)\log((1-r)/2)+\frac {n}
{2}(1+r)\log((1+r)/2)\\
-\frac {1} {2^n}\sum _{d=0} ^{\fl{n/2}}\frac {(n-2d+1)^2} {(n+1)}\binom {n+1}d
\log\la_d.
\end{multline}
The asymptotics of that expression can be determined in a similar way
to what was done for $0<r<1$. For the sake of brevity, we omit the derivation.
The result is (\ref{a4}). Actually, it is possible to rearrange the
computations that we did for $0<r<1$, so that in the limit $r\downarrow0$ they
give a proof of (\ref{a4}). This last observation justifies the
claim that the error term in (\ref{a3}) is uniform in $u$ and $r$,
$0\le r\le 1-\ep$ (i.e., including $r=0$), for any fixed positive
$\ep$.

\medskip
This completes the proof of the Theorem.
\quad \quad \qed
\medskip

Now, we list the summations which were used in the proof of the Theorem.

\begin{lemma}\label{la1}
We have the following summations:
\begin{equation}\label{a8}
\frac {1} {2^{n+1}r}\sum _{d=0} ^{n+1}\frac {(n-2d+1)}
{(n+1)} \binom {n+1}d
(1+r)^{n+1-d}(1-r)^d=1.
\end{equation}
\begin{equation}\label{a9}
\frac {1} {2^{n+1}r}\sum _{d=0} ^{n+1}\frac {(n-2d+1)}
{(n+1)} \binom {n+1}d
(1+r)^{n+1-d}(1-r)^dd=\frac {(1-r)(nr-1)} {2r}.
\end{equation}
\begin{equation}\label{a10}
\frac {1} {2^{n+1}r}\sum _{d=-1} ^{n+1}\frac {(n-2d+1)}
{(n+1)(n+2)} \binom {n+2}{d+1}
(1+r)^{n+1-d}(1-r)^d=\frac {2(1+2r+nr)} {(n+1)(n+2)r(1-r)}.
\end{equation}
\begin{multline}\label{a11}
 \frac {1} {2^{n+1}r}\sum _{d=0} ^{n+1}\frac {(n-2d+1)}
{(n+1)} \binom {n+1}d
(1+r)^{n+1-d}(1-r)^d  {(1/2-u+d)}\\
={\frac{-1 + 2 r + n r - n {r^2} - 2 r u} {2 r}}.
\end{multline}
\begin{multline}\label{a12}
 \frac {1} {2^{n+1}r}\sum _{d=0} ^{n+1}\frac {(n-2d+1)}
{(n+1)} \binom {n+1}d
(1+r)^{n+1-d}(1-r)^d \\
\hskip3cm\cdot {(1/2-u+d)}\(1+d-u-\frac {n(1-r)} {2}\)\\
={\frac{-5 - n + 7 r + 5 n r - 3 n {r^2} - n {r^3} + 4 u - 10 r u - 
     2 n r u + 2 n {r^2} u + 4 r {u^2}} {4 r}}.
\end{multline}
\begin{multline}\label{a13}
 \frac {1} {2^{n+1}r}\sum _{d=0} ^{n+1}\frac {(n-2d+1)}
{(n+1)} \binom {n+1}d
(1+r)^{n+1-d}(1-r)^d \\
\hskip3cm\cdot {(1/2-u+d)}\(1+d-u-\frac {n(1-r)} {2}\)^2\\
=\frac {1} {8r}(
-22 - 9 n + 26 r + 24 n r + {n^2} r - 5 n {r^2} - {n^2} {r^2} - 
     8 n {r^3} - {n^2} {r^3} - 2 n {r^4} + {n^2} {r^4} + 32 u \\
+ 
     4 n u - 48 r u - 22 n r u + 12 n {r^2} u + 6 n {r^3} u - 
     12 {u^2} + 32 r {u^2} + 4 n r {u^2} - 4 n {r^2} {u^2} - 
     8 r {u^3}).
\end{multline}
\begin{multline}\label{a14} 
\frac {1} {2^{n+1}r}\sum _{d=0} ^{n+1}\frac {(n-2d+1)}
{(n+1)} \binom {n+1}d
(1+r)^{n+1-d}(1-r)^d \\
\hskip3cm\cdot {(1/2-u+d)}\(1+d-u-\frac {n(1-r)} {2}\)^3\\
=\frac {1} {16r}(
-92 - 61 n - 3 {n^2} + 100 r + 105 n r + 15 {n^2} r + 19 n {r^2} - 
     4 {n^2} {r^2} - 35 n {r^3} - 20 {n^2} {r^3} - 22 n {r^4} \\+ 
     7 {n^2} {r^4} - 6 n {r^5} + 5 {n^2} {r^5} + 188 u + 60 n u - 
     228 r u - 162 n r u - 6 {n^2} r u + 20 n {r^2} u \\+ 
     6 {n^2} {r^2} u + 66 n {r^3} u + 6 {n^2} {r^3} u + 
     16 n {r^4} u - 6 {n^2} {r^4} u - 132 {u^2} - 12 n {u^2} + 
     204 r {u^2}\\
 + 72 n r {u^2} - 36 n {r^2} {u^2} - 
     24 n {r^3} {u^2} + 32 {u^3} - 88 r {u^3} - 8 n r {u^3} + 
     8 n {r^2} {u^3} + 16 r {u^4}).
\end{multline}

\end{lemma}

\medskip
{\sc Proof}. In all the cases, the sums can be split into several simpler
sums, each  of which can itself be summed using
 the binomial theorem.\quad \quad \qed
\medskip

\begin{lemma}\label{la2}
For fixed $r$ with $0<r<1$
we have the following asymptotic expansion:
\begin{multline}\label{a15} 
\frac {1} {2^{n+1}r}\sum _{d=0} ^{n+1}\frac {(n-2d+1)}
{(n+1)} \binom {n+1}d
(1+r)^{n+1-d}(1-r)^d (1/2-u+d)\log(1+d-u)\\
=\(\frac {n} {2}(1-r)+1-u-\frac {1} {2r}\)\big(\log
n+\log(1-r)-\log2\big)
\\+\frac {7} {4}-u+\frac {r} {4}-\frac {1} {2r}+O\(\frac {1} {n}\).
\end{multline}

\end{lemma}

\medskip
{\sc Proof}.
We start with the expansion
\begin{multline}\notag
\log(1+d-u)=\log\(\frac {n(1-r)} {2}\)+\log\(1+\frac {2}
{n(1-r)}\(1+d-u-\frac {n(1-r)} {2}\)\)\\
=\log n+\log(1-r)-\log 2+\frac {2}
{n(1-r)}\(1+d-u-\frac {n(1-r)} {2}\)\\
-\frac {2}
{n^2(1-r)^2}\(1+d-u-\frac {n(1-r)} {2}\)^2\\
+O\(\frac {1} 
{n^3(1-r)^3}\(1+d-u-\frac {n(1-r)} {2}\)^3\).
\end{multline}
(It is at this point that we must have $r<1$.)
If we use this expansion in the left-hand side of (\ref{a15}) and
subsequently use (\ref{a11})--(\ref{a14}) to evaluate the resulting
sums, we obtain 
\begin{multline}\label{a16}
\frac {1} {2^{n+1}r}\sum _{d=0} ^{n+1}\frac {(n-2d+1)}
{(n+1)} \binom {n+1}d
(1+r)^{n+1-d}(1-r)^d (1/2-u+d)\log(1+d-u)\\
=\(\frac {n} {2}(1-r)+1-u-\frac {1} {2r}\)\big(\log
n+\log(1-r)-\log2\big)
\\+\(2-u+\frac {r} {2}-\frac {1} {2r}\)-\frac {1+r} {4}+O\(\frac {1} {n}\).
\end{multline}
Simplifying easily, we obtain (\ref{a15}).\quad \quad \qed
\medskip

\subsection{Asymptotics of the von Neumann
entropies of the Bayesian density
matrices $\ze_n(u)$}\label{s2.5}
The main result of this section describes the asymptotics of the
von Neumann
entropy (\ref{eq:1}) of $\ze_n(u)$. In view of the explicit description of
the eigenvalues of $\ze_n(u)$ and their multiplicities in
Theorem~\ref{t2}, this entropy equals 
$$-\sum _{d=0} ^{\fl{n/2}}\frac {(n-2d+1)^2} {(n+1)}\binom {n+1}d
\la_d\log\la_d,$$
with $\la_d$ being given by (\ref{e14}).
\begin{theorem}\label{t11}
We have the following asymptotic expansion:
\begin{multline}\label{B2}
 -\sum _{d=0} ^{\fl{n/2}}\frac {(n-2d+1)^2} {(n+1)}\binom {n+1}d
\la_d \log\la_d=\\
 n \left( {\frac{-7 + 5 u} 
       {2 \left( 2 - u \right)  \left( 1 - u \right) }} + 
     {\psi}(5 - 2 u) - {\psi}(1 - u) \right)+
 \frac {3} {2}\log n+ 
  \left( -{\frac 7 2} + 2 u \right)  \log 2 \\-
{\frac{14 - 20 u + 7 {u^2}} 
    {2 \left( 2 - u \right)  \left( 1 - u \right) }}  + 
  \log \big({\Gamma}(1 - u)\big) - \log \big({\Gamma}({5/ 2} -
u)\big)\\
+(2-2u)(\psi(5-2u)-\psi(1-u))+O\(\frac {1} {n^{1-u}}\),
\end{multline}
where $\psi(x)$ is the digamma function,
$$\psi(x)=\frac {\frac {d} {dx}\Gamma(x)} {\Gamma(x)}.$$

\end{theorem}

\medskip
The proof of the Theorem depends on a few summations, which we now list.

\begin{lemma}\label{lB2}
We have the following summations:
\begin{equation}\label{B3}
\sum _{d=0} ^{n+1}\frac {(n-2d+1)^2} {(n+1)}\binom {n+1}d\la_d=2.
\end{equation}
\begin{equation}\label{B4}
\sum _{d=1} ^{n+1}{(n-2d+1)^2} \binom {n}{d-1}\la_d=n+1.
\end{equation}
\begin{equation}\label{B5}
\sum _{d=-1} ^{n+1}\frac {(n-2d+1)^2} {(n+1)(n+2)}\binom
{n+2}{d+1}\la_d=\frac {2(n+3)(2u-3)} {(n+1)(n+2)u}.
\end{equation}
\begin{multline}\label{B6}
\sum _{d=0} ^{n+1}\frac {(n-2d+1)^2} {(n+1)}\binom {n+1}d\\
\cdot
\frac {1} {2^n}\frac {\Ga(5/2-u)\,\Ga(2+n-d-u)\,\Ga(1+\al+d-u)}
{\Ga(5/2+n/2-u)\,\Ga(2+n/2-u)\,\Ga(1-u)}(d-u+1/2)\\
=( 48 + 64 {\alpha} + 25 {{{\alpha}}^2} + 3 {{{\alpha}}^3} + 40 n + 
       66 {\alpha} n + 37 {{{\alpha}}^2} n + 5 {{{\alpha}}^3} n + 
       8 {n^2} + 14 {\alpha} {n^2} + 8 {{{\alpha}}^2} {n^2}\\ + 
       2 {{{\alpha}}^3} {n^2} - 152 u - 138 {\alpha} u - 
       34 {{{\alpha}}^2} u - 2 {{{\alpha}}^3} u - 92 n u - 
       92 {\alpha} n u - 32 {{{\alpha}}^2} n u - 
       2 {{{\alpha}}^3} n u \\
- 12 {n^2} u - 10 {\alpha} {n^2} u - 
       2 {{{\alpha}}^2} {n^2} u + 176 {u^2} + 100 {\alpha} {u^2} + 
       12 {{{\alpha}}^2} {u^2} + 68 n {u^2} + 32 {\alpha} n {u^2}\\
 + 
       4 {{{\alpha}}^2} n {u^2} + 4 {n^2} {u^2} - 88 {u^3} - 
       24 {\alpha} {u^3} - 16 n {u^3} + 16 {u^4} )  \\
\times{\frac{
     \Gamma({ \textstyle 5 - 2 u})\,  
     \Gamma({ \textstyle 3 + {\alpha} + n - 2 u})\,  
     \Gamma({ \textstyle 1 + {\alpha} - u}) } 
   {4\, \Gamma({ \textstyle 5 + {\alpha} - 2 u})  \,
     \Gamma({ \textstyle 4 + n - 2 u})\,  \Gamma({ \textstyle 3 - u}) }}.
\end{multline}
\begin{multline}\label{B7} 
\sum _{d=0} ^{n+1}\frac {(n-2d+1)^2} {(n+1)}\binom {n+1}d\la_d
(d-u+1/2)\psi(1+d-u)\\
={\frac{32 + 33 n + 7 {n^2} - 69 u - 46 n u - 5 {n^2} u + 50 {u^2} + 
      16 n {u^2} - 12 {u^3}} 
    {2 \left( 2 - u \right)  \left( 1 - u \right)  
      \left( 3 + n - 2 u \right) }} \\
+ (n+2-2u)(\psi(1-u)+\psi(n+3-2u)-\psi(5-2u)).
\end{multline}

\end{lemma}
\medskip

{\sc Proof}. Identities (\ref{B3}), (\ref{B4}), (\ref{B5}), (\ref{B6}) 
are proved by splitting the
sums appropriately so that each part can be summed by means of Gau\ss'
$_2F_1$ summation. Identity (\ref{B7}) follows from (\ref{B6}) 
by differentiating
with respect to $\al$ and then setting $\al=0$.\quad \quad \qed
\medskip

{}From (\ref{B7}) we can deduce the following important estimation.
The result and its proof were kindly reported to us by Peter Grabner.
\begin{lemma}\label{lB3} 
We have the asymptotic expansion:
\begin{multline}\label{B8} 
\sum _{d=0} ^{n+1}\frac {(n-2d+1)^2} {(n+1)}\binom {n+1}d\la_d
(d-u+1/2)\log(1+d-u)\\
=     n \left( \log n + {\frac{7 - 5 u} 
         {2 \left( 2 - u \right)  \left( 1 - u \right) }} - 
       {\psi}(5 - 2 u) + {\psi}(1 - u) \right)+
\left( 2 - 2 u \right)\log n  \\ + 
    {\frac{26 - 46 u + 25 {u^2} - 4 {u^3}} 
      {2 \left( 2 - u \right)  \left( 1 - u \right) }} + 
    \left( -2 + 2 u \right)  {\psi}(5 - 2 u) + 
    \left( 2 - 2 u \right)  {\psi}(1 - u) 
  + 
  O\({\frac1 {n^{1-u}}}\).\\
\end{multline}
\end{lemma}
\medskip

{\sc Proof}. We use the asymptotic expansion
\begin{equation}\label{B9}
\psi(z)=\log(z)-\frac {1} {2z}+O\(\frac {1}
{z^2}\).
\end{equation}
In particular, this gives
$$\psi(1+d-u)=\log(1+d-u)-\frac {1} {2(d+1)}+O\(\frac {1}
{(d+1)(d+2)}\)$$
and
$$\psi(n+3-2u)=\log(n+2-2u)+\frac {1} {2(n+2-2u)}+O\(\frac {1} {n^2}\).$$
Using this in (7), we obtain
\begin{multline}\label{B10} 
\sum _{d=0} ^{n+1}\frac {(n-2d+1)^2} {(n+1)}\binom {n+1}d\la_d
(d-u+1/2)\log(1+d-u)\\
=\frac {1} {2}\sum _{d=0} ^{n+1}\frac {(n-2d+1)^2} {(n+1)}\binom {n+1}d\la_d
\frac {(d-u+1/2)} {d+1}\hskip2cm\\
+O\(
\sum _{d=0} ^{n+1}\frac {(n-2d+1)^2} {(n+1)}\binom {n+1}d\la_d
\frac {(d-u+1/2)} {(d+1)(d+2)}\)\\
+
\left( 2 - 2 u \right)\log n   - 
    {\frac{-22 + 40 u - 23 {u^2} + 4 {u^3}} 
      {2 \left( 2 - u \right)  \left( 1 - u \right) }} + 
    \left( -2 + 2 u \right)  {\psi}(5 - 2 u) + 
    \left( 2 - 2 u \right)  {\psi}(1 - u) \\
+     n \left( \log n + {\frac{7 - 5 u} 
         {2 \left( 2 - u \right)  \left( 1 - u \right) }} - 
       {\psi}(5 - 2 u) + {\psi}(1 - u) \right)  + 
  { O}\({\frac1 n}\)
\end{multline}
In the first expression on the right-hand side of (\ref{B10}) we use the
trivial identity
$$\frac {d-u+1/2} {d+1}=1-\frac {u+1/2} {d+1},$$
to split the expression into two sums, one of which can be evaluated
by means of (\ref{B3}). The other sum equals basically $-(u+1/2)$ times the
sum on the left-hand side of (\ref{B5}). What is missing is the summand for
$d=-1$. By (\ref{B5}), the complete sum is of the order $O(1/n)$. Using
Stirling's formula it is seen that the summand for $d=-1$ is of the
order $O(1/n^{1-u})$. So, combining everything, the first expression in
(\ref{B10}) equals $1+O(1/n)+O(1/n^{1-u})=1+O(1/n^{1-u})$.
For the second expression, we do a similar partial
fraction expansion in order to apply (\ref{B5}). The result is that
this second expression is of the order $O(1/n^{1-u})$.
This establishes the Lemma.\quad \quad \qed
\medskip

Now we are in the position to prove the Theorem.
\medskip

{\sc Proof of the Theorem}.
Since $\la_{n+1-d}=\la_d$, an equivalent expression for the
left-hand side in (\ref{B2}) is
\begin{equation}\label{B11}
-\frac {1} {2}\sum _{d=0} ^{n+1}\frac {(n-2d+1)^2} {(n+1)}\binom {n+1}d
\la_d \log\la_d.
\end{equation}
Now, we expand the logarithm according to the addition rule to obtain
\begin{align}\notag 
-\frac {1} {2}&\sum _{d=0} ^{n+1}\frac {(n-2d+1)^2}
{(n+1)}\binom {n+1}d \la_d\\
\notag
&\hskip3cm\cdot\log\(\frac {1} {2^n}
\frac {\Ga(5/2-u)}
{\Ga(5/2+n/2-u)\,\Ga(2+n/2-u)\,\Ga(1-u)}\)\\
\notag
&\hskip.8cm -\frac {1} {2}\sum _{d=0} ^{n+1}
{\frac {(n-2d+1)^2} {(n+1)}\binom
{n+1}d}\la_d\big(\log\Ga(1+d-u)+\log\Ga(2+n-d-u)\big).
\end{align}

The first sum in this expression can be evaluated by means of (\ref{B3}).
Therefore, we obtain for the expression on the left-hand side of (\ref{B2})
\begin{multline}\label{B13}
n\log2-\log\Ga(5/2-u)+\log\Ga(1-u)
+\log\Ga(5/2+n/2-u)\\
 +\log\Ga(2+n/2-u)-\sum _{d=0} ^{n+1}
{\frac {(n-2d+1)^2} {(n+1)}\binom
{n+1}d}\la_d\log\Ga(1+d-u).
\end{multline}

The only difficulty in obtaining the asymptotics of expression
(\ref{B13}) stems from the sum. In this sum, 
we use Stirling's formula
$$\log\Ga(x)=(x-1/2)\log x-x+\frac {1} {2}\log 2+\frac {1}
{2}\log\pi+O\(\frac {1} {x}\)$$
to get
\begin{align}\notag 
&\sum _{d=0} ^{n+1}
{\frac {(n-2d+1)^2} {(n+1)}\binom
{n+1}d}\la_d\log\Gamma(1+d-u)\\
\notag
&=\sum _{d=0} ^{n+1}{\frac {(n-2d+1)^2} {(n+1)}\binom
{n+1}d}\la_d\\
\notag
&\hskip1cm\cdot\Big((1/2+d-u)\log(1+d-u)-1+u-d+\frac {1} {2}\log 2+\frac {1}
{2}\log\pi+O\Big(\frac {1} {d+1}\Big)\Big)\\
\notag
&=\sum _{d=0} ^{n+1}
{\frac {(n-2d+1)^2} {(n+1)}\binom
{n+1}d}\la_d(u-1+\frac {1} {2}\log 2+\frac {1} {2}\log\pi) \\
\notag&\hskip2cm-\sum _{d=0} ^{n+1}
{ {(n-2d+1)^2} \binom {n}{d-1}}
 +O\(\sum _{d=0} ^{n+1}
{\frac {(n-2d+1)^2} {(n+1)(n+2)}\binom
{n+2}{d+1}}\)\\
\label{B14}
&\hskip2cm
+\sum _{d=0} ^{n+1}{\frac {(n-2d+1)^2} {(n+1)}\binom
{n+1}d}(1/2-u+d)\log(1-u+d).
\end{align}

The first expression on the right-hand side of (\ref{B14}) simplifies by means
of (\ref{B3}), the second by means of (\ref{B4}). For the $O(.)$ term
we use (\ref{B5}). In
fact, the sum on the left-hand side of (\ref{B5}) differs from the sum in the
$O(.)$ term only by the summand for $d=-1$. This summand is of the
order $O(1/n^{1-u})$, as is seen by Stirling's formula. Putting
everything together, we obtain 
\begin{multline}\label{B15} 
\sum _{d=0} ^{n+1}
{\frac {(n-2d+1)^2} {(n+1)}\binom
{n+1}d}\la_d\log\Gamma(1+d-u)\\
=2u-2+\log2+\log\pi-(n+1)+O\(\frac {1} {n^{1-u}}\)\\
+\sum _{d=0} ^{n+1}{\frac {(n-2d+1)^2} {(n+1)}\binom
{n+1}d}\la_d(1/2-u+d)\log(1-u+d).
\end{multline}
When we use this in (\ref{B13}),
 apply Lemma~3 to the remaining sum, and
simplify, we finally arrive at (\ref{B2}).\quad \quad \qed

\section{Comparison of our asymptotic redundancies
for the one-parameter family $q(u)$ with those of Clarke
and Barron}\label{s3}
Let us, first, compare the formula (\ref{eq:4}) for the asymptotic redundancy
of Clarke and Barron to that derived here (\ref{a3}) for the two-level
quantum systems, in terms of the one-parameter family of probability
densities $q(u)$, $ -\infty < u < 1$, given in (\ref{eq:10}).
Since the unit ball or Bloch sphere of such systems is three-dimensional
in nature, we are led to set the dimension $d$ of the parameter space in
(\ref{eq:4}) to 3. 
The quantum Fisher information matrix $I(\theta)$ 
for that case was taken to be (\ref{eq:8}), while the role
of the probability function $w(\theta)$ is played by $q(u)$.
Under these substitutions, it was seen in the Introduction that 
formula (\ref{eq:4}) reduces to (\ref{eq:12}).
Then, we see that for 
$0 \le r < 1$, the formulas (\ref{a3}) and
(\ref{eq:12})
coincide except for the presence of the
monotonically decreasing (nonclassical/quantum) term
$\frac1  {2r} \log \(\frac{1-r} { 1+r}\)$ (see Figure~2 
for a plot of this term --- $\log 2 \approx .693147$ ``nats'' of
information equalling one ``bit'')
in (\ref{a3}). (This term would
 have to be replaced by $-1$ --- that is,
its limit for $r \rightarrow 0$ --- to give (\ref{eq:12}).)
In particular, the order of magnitude, ${\frac3  2} \log n$,
is precisely the same in both formulas.
For the particular case $r=0$, the asymptotic 
formula (\ref{a3}) (see (\ref{a4}))
precisely coincides with (\ref{eq:12}).

\vskip10pt
\vbox{
\centerline{
\hbox{\psfig{file=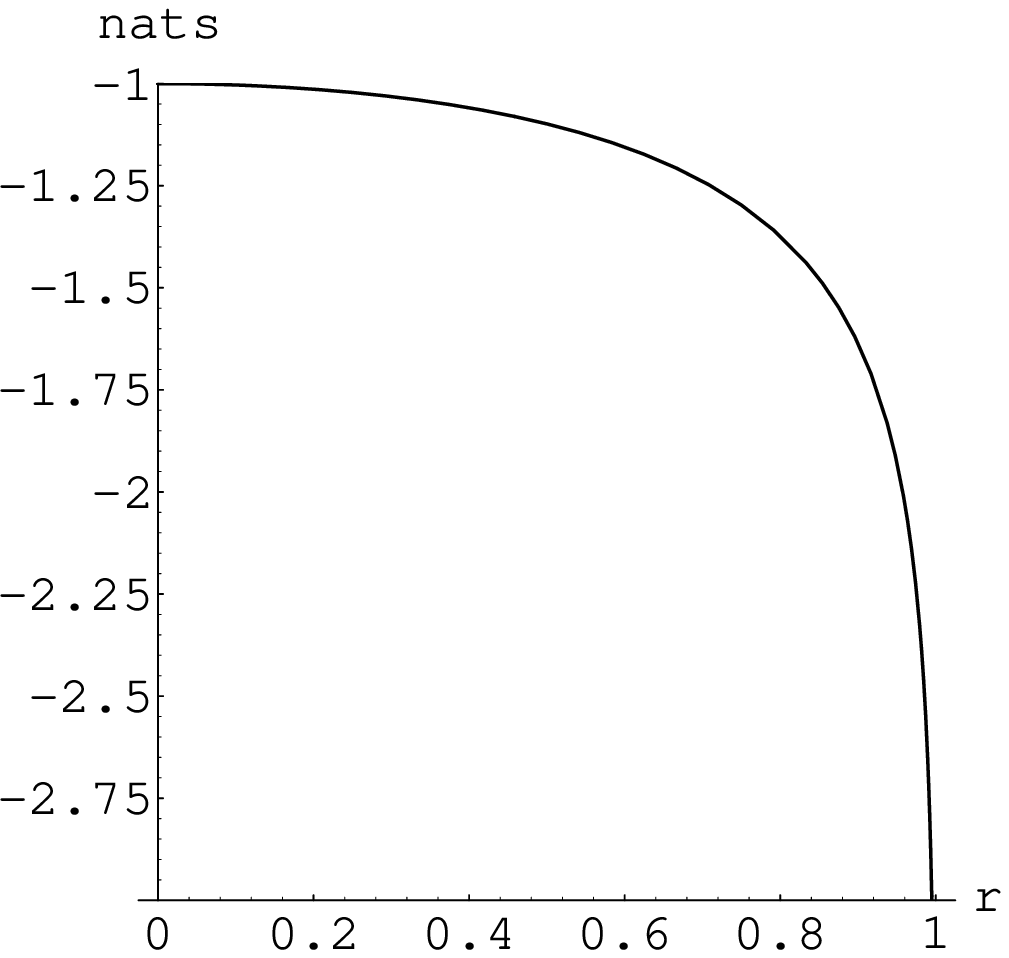,height=7cm}}
}
\centerline{\small Nonclassical/quantum term ($\frac1  {2 r} \log
\frac{1-r} {1+r}$) 
in the quantum asymptotic redundancy (\ref{a3})}
\vskip5pt
\centerline{\small Figure 2}
}
\vskip10pt

In the case $r=1$, however, i.e., when we consider the
boundary of the parameter space (represented by the unit sphere), 
the situation is slightly tricky.
Due to the fact that the formula of Clarke and Barron 
holds only for interior points of the parameter space, we cannot expect
that, in general, our formula will resemble that of Clarke and Barron.
However, if the probability density, $q(u)$, is concentrated on
the boundary of the sphere, then we may disregard the interior of the
sphere, and may consider the boundary of the sphere as the {\em true}
parameter space. This parameter space is {\em two-dimensional\/} and
consists of interior points throughout.
Indeed, the probability density $q(u)$ is concentrated on the
boundary of the sphere if we choose $u=1$ since, as
we remarked in the Introduction, in the limit $u\to1$, the
distribution determined by $q(u)$ tends to the uniform distribution
over the boundary of the sphere. Let us, again, (naively) attempt to apply
Clarke and Barron's formula (\ref{eq:4}) to that case. We parameterize
the boundary of the sphere by polar coordinates 
$(\vartheta,\phi)$,
\begin{align}\notag x&=\sin\vartheta \cos\ph\\
y&=\sin\vartheta \sin\ph\notag\\
z&=\cos\vartheta,\notag\\
0\le \ph\le {}&2\pi,\ 0\le\vartheta\le\pi.\notag
\end{align}
The probability density induced by $q(u)$ in the limit $u\to1$ then
is $\sin\vartheta/4\pi$, the density of the uniform distribution.
Using
\cite[eq.~8]{fuj}, the quantum (symmetric logarithmic derivative)
Fisher information matrix turns out to be
\begin{equation}
\begin{pmatrix}1&0\\
0&\sin^2\vartheta\end{pmatrix}\quad,
\end{equation}
its determinant equalling, therefore, $\sin^2\vartheta$. So, 
setting $d=2$ and substituting
$\sin\vartheta/4\pi$ for $w(\th)$ and $\sin^2\vartheta$ for $I(\th)$ in
(\ref{eq:4})
gives $\log n+\log 2-1$. On the other hand, our formula (\ref{a5}),
for $u=1$, gives $\log n$. So, again, the terms differ only by a
constant. In particular, the order of magnitude is again the same.

Let us now focus our attention on the asymptotic minimax
redundancy (\ref{eq:3}) of Clarke and Barron. If in (\ref{eq:3}) we again set
$d$ to 3, we obtain (\ref{eq:11}).
Clarke and Barron prove that this minimax expression is only attained by the
(classical) Jeffreys' prior.
In order to derive its quantum counterpart
--- at least, a restricted (to the family $q(u)$) version --- we
 have to determine the behavior of
\begin{equation}\label{d5}
\min_{-\infty<u<1}\max_{0\le r\le 1}S(\overset n\otimes \rho, \ze_n(u))
\end{equation}
for $n \rightarrow \infty$.
We are unable to proceed in a fully rigorous manner. However,
from computational data we conjecture that
\begin{equation}\label{d6}
\max_{0\le r\le 1}S(\overset n\otimes \rho, \ze_n(u))
\end{equation}
is always attained at $r=0$
(corresponding to the fully mixed state)
 or $r=1$ (corresponding to a pure state). Assuming the validity
of this conjecture,  the maximum $u_n$ in (\ref{d6}) is 
a value for which $S( \overset n\otimes \rho,\zeta_{n}(u))|_{r=0}$
equals $S(\overset n\otimes \rho,\zeta_{n}(u))|_{r=1}$.
Then we are able to prove that $\lim_{n\to \infty} u_{n} = .5$. 

Namely, by our assumption we have
\begin{equation}\label{d6a}
S( \overset n\otimes \rho,\zeta_{n}(u_{n}))|_{r=0}=
S(\overset n\otimes \rho,\zeta_{n}(u_{n}))|_{r=1},
\end{equation}
for any $n$.
Let $(u_{n_k})_{k=1,2,\dots}$ be a subsequence of
the sequence $(u_n)$ which converges to some $u_0$, $-\infty\le u_0\le 1$.
Note that we allow $u_0=-\infty$ and $u_0=1$. Therefore, there is
always such a subsequence. Because of (\ref{d6a}) we must have
\begin{equation}\label{d6b}
\lim_{k\to\infty}\frac 
{S( \overset {n_k}\otimes \rho,\zeta_{n_k}(u_{n_k}))|_{r=0}} {\log
n_k}=
\lim_{k\to\infty}\frac 
{S(\overset {n_k}\otimes \rho,\zeta_{n_k}(u_{n_k}))|_{r=1}} {\log
n_k}.
\end{equation}
By (\ref{a4}) and the fact that the error term in (\ref{a4}) is
uniform in $u$, we know that the left-hand side in (\ref{d6b}) is $3/2$. On
the other hand, by (\ref{a5}) and 
the fact that the error term in (\ref{a5}) is
uniform in $u$,
the right-hand side in (\ref{d6b}) equals
$\lim_{k\to\infty}(2-u_{n_k})$. Hence, we must have
$\lim_{k\to\infty}u_{n_k}=.5$. Thus, every convergent subsequence of
$(u_n)$ (including those which converge to $-\infty$ or $1$, the
boundary points of the interval of possible values of $u_n$)
converges to $.5$. Hence, the complete sequence $(u_n)$ converges to
$.5$, establishing our claim.
Since we have regarded $q(.5)$, that is (\ref{eq:9}), as
the quantum counterpart of the
 Jeffreys' prior (because, by analogy with the classical situation, it is
the normalized square root of the determinant of the
quantum Fisher information matrix, $\sqrt {{\det I(\theta)}}$),
this result could be considered to be
 fully parallel to that of Clarke and Barron.

We now concern ourselves with the asymptotic {\it maximin} redundancy.
Clarke and Barron \cite{cl1,cl2} prove that the maximin redundancy is
attained asymptotically, again, by the Jeffreys' prior.
To derive the quantum counterpart of the maximin redundancy within
our analytical framework, we would have to calculate
\begin{equation}\label{d7}
\max _w\min_{Q_n}\int_{x^2+y^2+z^2\le 1}S(\overset
n\otimes\rho,Q_n)\,w(x,y,z)\,dx\,dy\,dz,
\end{equation}
where $Q_{n}$ varies over the $(2^{2 n}-1)$-dimensional convex set
of $2^n \times 2^n$ density matrices and $w$ varies over all
probability 
densities over the unit ball.
In the classical case, due to a result of Aitchison 
\cite[pp.~549/550]{ait},
the minimum is achieved by setting
$Q_{n}$ to be the Bayes estimator,
i.e., the average of all possible $Q_n$'s with respect to the given
probablity distribution.
In the quantum domain the same 
assertion is true. For the sake of completeness, we include the proof
in the Appendix.
We can, thus, take the quantum
analog of the Bayes estimator 
to be the Bayesian density matrix $\zeta_{n}(u)$. 
That is, we set $Q_n=\ze_n(u)$ in (\ref{d7}).
Let us, for the moment, restrict the possible $w$'s over which
the maximum is to be taken to the family
$q(u)$, $-\infty < u < 1$. Thus, we consider
\begin{equation} \label{d7a}
\max _u\int_{x^2+y^2+z^2\le 1} S(\overset
n\otimes\rho,\ze_n(u))\,q(u)\,dx\,dy\,dz. 
\end{equation}
By the definition (\ref{eq:5}) of relative entropy, we have
\begin{multline}\notag
S(\overset
n\otimes\rho,\ze_n(u))=\Tr\(\overset n\otimes\rho\log \overset
n\otimes\rho\) -\Tr\(\overset n\otimes\rho\log \ze_n(u)\)\\
=n{\frac{(1-r)}  2}{ \log {\frac{(1-r)}  2}} +n{\frac {(1+r) }
2}{ \log {\frac{(1+r)}   2}}-\Tr\(\overset n\otimes\rho\log
\ze_n(u)\),
\end{multline}
the second line being due to (\ref{eq:7}). Therefore, we get
\begin{multline} \label{d7b}
\int_{x^2+y^2+z^2\le 1} S(\overset
n\otimes\rho,\ze_n(u))\,q(u)\,dx\,dy\,dz\\
=\(n\int_0^1\int_0^{\pi}\int_0^{2\pi}\(
{\frac{(1-r)}  2}{ \log {\frac{(1-r)}  2}} +{\frac {(1+r) }
2}{ \log {\frac{(1+r)}   2}}\)\)\,d\ph\,d\vartheta\,dr\\
-\Tr\big(\ze_n(u)\log\ze_n(u)\big)\\
=- n \left( {\frac{-7 + 5 u} 
       {2 \left( 2 - u \right)  \left( 1 - u \right) }} + 
     {\psi}(5 - 2 u) - {\psi}(1 - u) \right)
+S(\ze_n(u)).
\end{multline}
{}From Theorem~\ref{t11}, we know the asymptotics of the von Neumann
entropy $S(\ze_n(u))$. Hence, we find that the
expression (\ref{d7b}) is asymptotically
equal to
\begin{multline} \label{d8}
 \frac {3} {2}\log n+ 
  \left( -{\frac 7 2} + 2 u \right)  \log 2 \\-
{\frac{14 - 20 u + 7 {u^2}} 
    {2 \left( 2 - u \right)  \left( 1 - u \right) }}  + 
  \log \big({\Gamma}(1 - u)\big) - \log \big({\Gamma}({5/ 2} -
u)\big)\\
+(2-2u)(\psi(5-2u)-\psi(1-u))+O\(\frac {1} {n^{1-u}}\).
\end{multline}
We have to, first, perform the maximization required in (\ref{d7a}),
and then determine the asymptotics of the result. 
Due to the form of the asymptotics in (\ref{d8}), we can, in fact,
derive the proper result by proceeding in the reverse order. That is,
we first determine the asymptotics of $\int S(\overset
n\otimes\rho,\ze_n(u))\,q(u)\,dx\,dy\,dz$, which we did in (\ref{d8}),
and then we
maximize the $u$-dependent part in (\ref{d8}) with respect to $u$
(ignoring the error term).
(In Figure~3 we display this $u$-dependent part over the range
$[-0.2,1]$.)
Of course, we do the latter step by equating the first derivative of the
$u$-dependent part in (\ref{d8})
with respect to $u$ to zero and solving for $u$.
It turns out that this equation takes the appealingly simple form
\begin{equation} \label{maximin}
2(1-u)^3\big(\psi'(1-u)-\psi'(5/2-u)\big)=1.
\end{equation}
Numerically, we find  this equation to have the solution $u \approx .531267$,
at which the asymptotic maximin redundancy 
assumes the value ${\frac 3  2} \log n 
-1.77185+
O(1/n^{.468733})$. For $u=.5$, on the other hand, we have
for the asymptotic minimax redundancy,
${\frac 3  2} \log n
-2-\frac {1} {2}\log 2+\frac {1}
{2}\log\pi+
O(1/\sqrt n)=\frac {3} {2}\log n
-1.77421
+O(1/\sqrt n)$.
We must, therefore, conclude that --- in contrast to the
classical case \cite{cl1,cl2} --- our trial candidate ($q(.5)$) for
the quantum counterpart of Jeffreys' prior can not
serve as a ``reference prior,'' in the sense introduced
by Bernardo \cite{berna,bern}.

\vskip10pt
\vbox{
\centerline{
\psfig{file=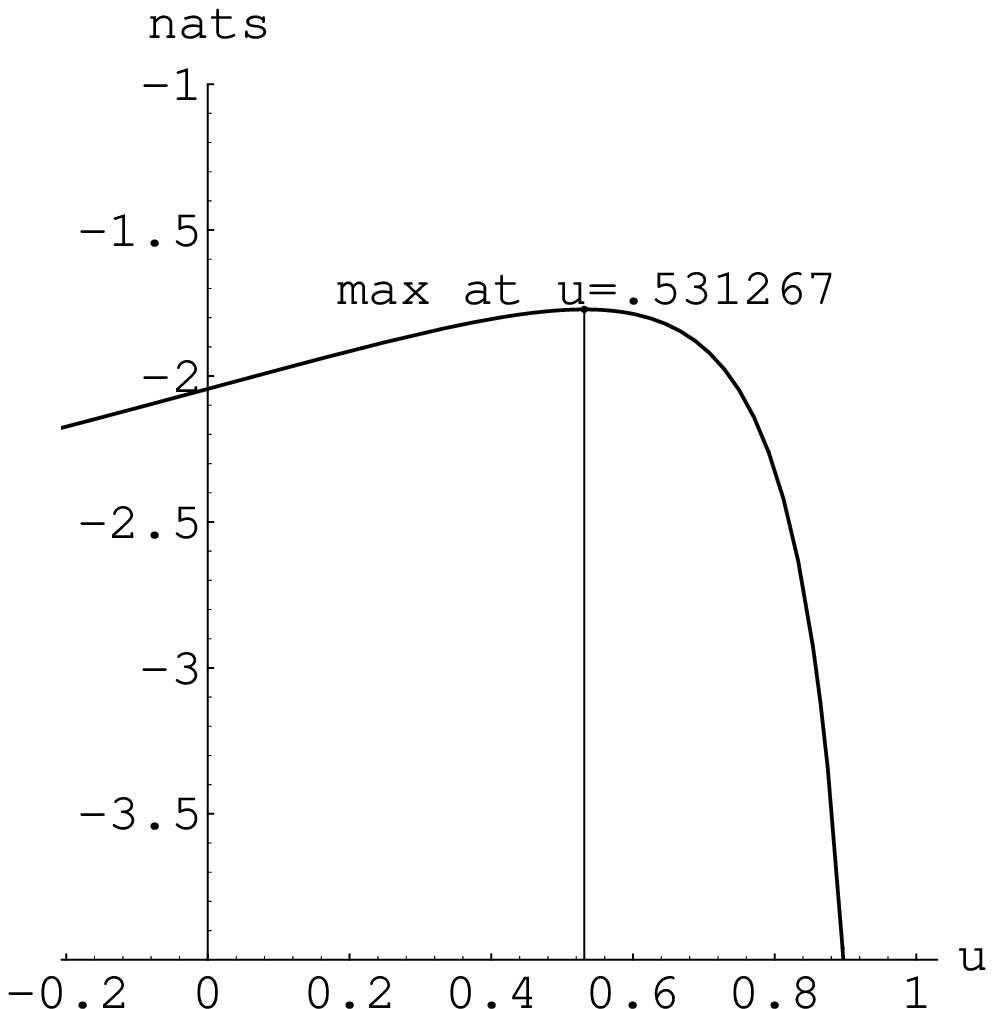,height=7cm}
}
\centerline{\small $u$-dependent part of the asymptotic Bayes
redundancy (\ref{d8})}
\vskip5pt
\centerline{\small Figure 3}
}
\vskip10pt

Since they are mixtures of product states, the matrices
$\zeta_{n}(u)$ are classically --- as opposed to EPR,
Einstein--Podolsky--Rosen --- correlated
\cite{werner}. Therefore, $S(\zeta_{n}(u))$ must not be less than
the sum of the von Neumann entropies 
of any set of reduced
density matrices obtained from it, through computation of partial traces.
For positive integers, $n_{1} + n_{2}+ \cdots = n$, the corresponding
reduced density matrices are simply $\zeta_{n_{1}(u)}, \zeta_{n_{2}(u)},
 \dots$,
due to the mixing
\cite[exercise~7.10]{belt}.
 Using these reduced density matrices, one can compute
{\it conditional} density matrices and quantum entropies \cite{cerf}.
Clarke and Barron \cite[p.~40]{cl1} have an alternative expression for
the redundancy in terms of conditional entropies, and it would be of
interest to ascertain whether a quantum analogue of this expression exists.

Let us note that the theorem of Clarke and Barron utilized the
uniform convergence property of the asymptotic expansion of the
relative entropy (Kullback-Leibler divergence). Condition 2 in their paper \cite{cl1} is, therefore,
crucial. It assumes --- as is typically the
 case classically --- that the matrix of
  second derivatives, $J(\theta)$, of the
relative entropy is identical to the Fisher information matrix $I(\theta)$.
In the quantum domain, however, in general, $J(\theta) \geq I(\theta)$,
where $J(\theta)$ is the matrix of second derivatives of the quantum relative 
entropy (\ref{eq:5}) and $I(\theta)$ is the symmetric logarithmic
derivative Fisher information matrix \cite{pe1,pe2}.
The equality holds only for special cases. For instance,
$J(\theta) > I(\theta)$ does hold if $r \neq 0$ for the situation
considered in this paper.
The volume element of the Kubo-Mori/Bogoliubov (monotone) metric \cite{pe1,pe2}
is given by $\sqrt {\det J(\theta)}$. This can be normalized
for the two-level quantum systems to be a member
($u=1/2$) of a one-parameter family of probability 
densities
\begin{equation} \label{eq:Kubo}
\frac{(1-u) \,\Gamma(5/2-u) \,r\, \log{\big((1+r)/(1-r)\big)}\, \sin \vartheta}
{\pi^{3/2}\, (3 - 2 u)\, \Gamma(1-u)\, (1-r^2)^u},\quad -\infty<u<1,
\end{equation}
and similarly studied, it is presumed,
in the manner of the family $q(u)$ (cf. (\ref{eq:10}) and (\ref{e7})) analyzed here. These
two families can be seen to differ --- up to the normalization
factor --- by the replacement of $\log {\big((1+r)/(1-r)\big)}$ in (\ref{eq:Kubo})
by, simply, $r$.
(These two last expressions are, of course, equal for $r=0$.)
In general, the volume element of a monotone metric over the two-level
quantum systems is of the form \cite[eq.~3.17]{pe1}
\begin{equation} \label{monotone}
\frac{r^2 \sin{\vartheta}}{f\big((1-r)/(1+r)\big) (1-r^2)^{1/2} (1+r)},
\end{equation}
where $f:\mathbb R^+ \rightarrow \mathbb R^+$ is an operator monotone function
such that $f(1) = 1$ and $f(t) = tf(1/t)$.
For $f(t)=(1+t)/2$, one recovers the volume element
($\sqrt {\det {I(\theta)}}$) of the
metric of the symmetric
logarithmic derivative, and for $f(t) = 
(t-1)/ {\log t}$, that
($\sqrt {\det {J(\theta)}}$) of the
Kubo-Mori/Bogoliubov metric \cite{pe3,pe1,pe2}.
(It would appear, then, that the only member of the family $q(u)$
proportional to a monotone metric is $q(.5)$, that is (\ref{eq:9}).
 The maximin result
we have obtained above corresponding to $u \approx .531267$ --- the solution
of (\ref{maximin}) --- would
 appear unlikely, then, to extend globally beyond the family.)
While $J(\theta)$ can be generated from the relative entropy (\ref{eq:5})
(which is a limiting case of the $\alpha$-entropies \cite{pe4}),
$I(\theta)$ is similarly obtained from \cite[eq.~3.16]{pe3}
\begin{equation} \label{jan}
 \Tr \rho_{1} (\log \rho_{1} - \log \rho_{2})^2 .
\end{equation}
It might prove of interest to repeat the general line of
analysis carried out in this paper, but with the use
of (\ref{jan}) rather than (\ref{eq:5}).
Also of importance might be an analysis in which the relative
entropy (\ref{eq:5}) is retained, but the family (\ref{eq:Kubo}) based
on the Kubo-Mori/Bogoliubov metric is used instead of $q(u)$.
Let us also indicate that if one equates the asymptotic
redundancy formula of Clarke and Barron (\ref{eq:4})
(using $w(\theta) = q(u)$) to that derived here
(\ref{a3}),
 neglecting the residual terms, solves for
$\det(I(\theta))$, and takes the square root of the result,
one obtains a prior of the form (\ref{monotone})
 based on the
monotone function $t^\frac{t}{1+t}$.

As we said in the Introduction, ideally we would like to start with a
(suitable well-behaved) {\it arbitrary\/} probability density
on the unit ball, determine the relative entropy of 
$\overset n\otimes \rho$ with respect to the average of 
$\overset n\otimes \rho$ over the probability density,
then find its asymptotics, and finally, among all such probability
densities, find the one(s) for which the minimax and maximin are
attained. In this regard, we wish to mention that a suitable
combination of results and computations from Sec.~\ref{s2} with basic
facts from representation theory of $SU(2)$ 
(cf\@. \cite{vilklim,biedlouck} for more information on that topic) 
yields the following result.

\begin{theorem}\label{t15}
Let $w$ be a spherically symmetric probability density on the unit
ball, i.e., $w=w(x,y,z)$ depends only on $r=\sqrt{x^2+y^2+z^2}$.
Furthermore, let $\hat \ze_n(w)$ be the average 
$\int_{x^2+y^2+z^2\le 1} \big(\overset n  \otimes \rho \big)\,w \,dx\,dy\,dz$.
Then the eigenvalues of 
$\hat\ze_n(w)$ are
\begin{equation}\label{e50}
\la_{d}=\frac {\pi} {2^{n-1}(n-2d+1)}\int_{-1}^{1}
r(1+r)^{n-d+1}(1-r)^d w(\v{r})\,dr, 
\quad d=0,1,\dots,\fl{\frac {n} {2}},
\end{equation}
with respective multiplicities
\begin{equation}\label{e51}
\frac {n-2d+1} {n+1}\binom {n+1}d,
\end{equation}
and corresponding eigenspaces as given by {\em(\ref{e18})}.

The relative entropy of $\overset n  \otimes \rho$ with respect to
$\hat\ze_n(w)$ is given by {\em(\ref{e32})}, with $\la_d$ as given in
{\em(\ref{e50})}.
\end{theorem}
We hope that this Theorem enables us to determine the
asymptotics of the relative entropy and, eventually, to 
find, at least within the family of spherically
symmetric probability densities on the unit ball, the corresponding
minimax and maximin redundancies.

\section{Summary}\label{s4}
Clarke and Barron \cite{cl1,cl2} (cf\@. \cite{ri}) have derived several
forms of asymptotic redundancy for arbitrarily parameterized
families of probability distributions.
We have been motivated to undertake this study by the possibility
that their results may generalize, in some yet not fully understood
fashion, to the quantum domain of noncommutative probability.
(Thus, rather than probability 
densities, we have been concerned
here with density matrices.)
We have only, so far, been able to examine this possibility in a
somewhat restricted manner.
By this, we mean that we have limited our consideration to two-level
quantum systems (rather than $n$-level ones, $n \geq 2$), and for the case
$n=2$, we have studied (what has proven
to be) an analytically tractable one-parameter family of possible prior
probability 
densities, $q(u)$, $-\infty <u<1$
(rather than the totality of arbitrary probability 
densities).
Consequently, our results can not be as definitive in
nature as
those of Clarke and Barron.
Nevertheless, the analyses presented here indicate
  that our trial candidate ($q(.5)$, that is (\ref{eq:9})) for the quantum counterpart
of the Jeffreys' prior plays a somewhat similarly
 privileged --- but less pronounced --- role
as in the classical case.

                                              Future research
might be devoted to expanding the family of probability distributions used
to generate the Bayesian density matrices for $n=2$, as well as
similarly studying the $n$-level quantum systems ($n>2$).
(In this
regard, we have examined the situation in which $n=2^m$, and the
only $n \times n$
density matrices considered are simply the tensor products of $m$ identical
$2 \times 2$ density matrices. Surprisingly,
for $m=2,3$, the associated trivariate
candidate  quantum Jeffreys' prior,
taken, as throughout this study,
 to be proportional to the volume elements of the metrics
of the symmetric logarithmic derivative (cf\@. \cite{sl4}),
have been found to be
   {\it improper} (nonnormalizable) over the
Bloch sphere. The minimality of such metrics is guaranteed, however,
only if ``the whole state space of a spin is parameterized'' \cite{pe1}.)
In all such cases, it will be of interest to evaluate the characteristics
of the relevant candidate quantum Jeffreys' prior {\it vis-\`a-vis} all other
members of the family
of probability distributions employed over the
$(n^2-1)$-dimensional
convex set of
$n \times n$ density matrices.

We have also conducted analyses parallel to those
reported above, but having,
{\it ab initio},
 set either $x$ or $y$ to zero in
the $2 \times 2$ density matrices (\ref{eq:6}). This,
then, places us in the realm
of real --- as opposed to complex ( standard or conventional)
quantum mechanics. (Of course, setting {\it both} $x$ and $y$
to zero would return us to a strictly classical situation, in which
the results of Clarke and Barron \cite{cl1,cl2}, as applied to binomial
distributions, would be directly applicable.)
Though we have --- on the basis
of detailed computations --- developed strong conjectures as to the nature
of the 
associated results, we have not,
at this stage of our investigation, yet succeeded in formally
demonstrating their validity.

In conclusion, again in analogy to classical results,
we would like to raise the possibility that the quantum
asymptotic redundancies derived here might prove of
value in deriving formulas for the {\it stochastic
complexity} \cite{ri,ri2} 
(cf.~\cite{svozil}) --- the shortest
description length --- of a string of $n$
{\it quantum} bits. The competing possible models for the data string
might be taken to be the $2 \times 2$ density matrices ($\rho$)
corresponding to different values of $r$, or equivalently,
different values of the von Neumann entropy, $S(\rho)$.

\section*{Appendix: The quantum Bayes estimator achieves the minimum
average entropy} 
Let $P_\th$, $\th\in\Th$, be a family of density matrices, and let
$w(\th)$, $\th\in\Th$, be a family of probability distributions. 

\begin{theorem}
The minimum
$$\min_Q\int w(\th)S(P_\th,Q)\,d\th,$$
taken over all density matrices $Q$, is achieved by $m=\int w(\th)P_\th\,d\th$.
\end{theorem}
{\sc Proof}. We look at the difference
$$\int w(\th)S(P_\th,Q)\,d\th-\int w(\th)S(P_\th,m)\,d\th,$$
and show that it is nonnegative. Indeed,
\begin{align*}
\int w&(\th)S(P_\th,Q)\,d\th-\int w(\th)S(P_\th,m)\,d\th\\&=
\int w(\th)\Tr(P_\th\log P_\th-P_\th\log Q)\,d\th-
\int w(\th)\Tr(P_\th\log P_\th-P_\th\log m)\,d\th\\
&=\int w(\th)\Tr\big(P_\th(\log m -\log Q)\big)\,d\th\\
&=\Tr\Big(\Big(\int w(\th)P_\th\,d\th\Big)(\log m -\log Q)\Big)\\
&=Tr\big(m(\log m -\log Q)\big)\\
&=S(m,Q)\ge0,
\end{align*}
since relative entropies are nonnegative \cite{oh}. \quad \quad \qed

\medskip

\section*{Acknowledgments}
Christian Krattenthaler did part of this research at the 
Mathematical Sciences Research Institute, Berkeley, during the
Combinatorics Program 1996/97.
Paul Slater would like to express appreciation
to the Institute for Theoretical Physics for
computational support. This research was undertaken, in part, to respond to
concerns (regarding the rationale for the
presumed quantum Jeffreys' prior)
 conveyed to him by Walter Kohn
and members of the informal seminar group he leads.
The co-authors are grateful 
to Helmut Prodinger and Peter Grabner for their hints regarding 
the asymptotic computations, to Ira Gessel for
bringing them into initial contact {\it via} the Internet, and
to A. R. Bishop and an anonymous referee of \cite{sl4}.

\end{document}